\documentclass{amsart}

\usepackage{amsmath, amsthm, amssymb, booktabs, multirow, graphicx,
  booktabs, bm, float}

\usepackage{mathtools}
\usepackage{enumitem}
\usepackage{dialogue}
\usepackage{bbold} 
\usepackage[stable]{footmisc}
\usepackage{tikz,nicefrac}
\usetikzlibrary{calc,matrix,arrows,decorations.markings}
\usepackage[colorlinks]{hyperref}
\usetikzlibrary{
  graphs,
  graphs.standard
}

\newtheorem{defin}{Definition}[section]

\newtheorem{theorem}[defin]{Theorem}

\newcommand{\C}{\mathbb{C}}
\newcommand{\R}{\mathbb{R}}

\newcommand{\Kcal}{\mathcal{K}}

\newcommand{\Ccal}{\mathcal{C}}

\newcommand{\defi}[1]{\textit{#1}}

\makeatletter
\newcommand{\bigperp}{%
  \mathop{\mathpalette\bigp@rp\relax}%
  \displaylimits
}

\newcommand{\bigp@rp}[2]{%
  \vcenter{
    \m@th\hbox{\scalebox{\ifx#1\displaystyle2.1\else1.5\fi}{$#1\perp$}}
  }%
}
\makeatother

\DeclareMathOperator{\vol}{vol}

\DeclareMathOperator{\cone}{cone}

\DeclareMathOperator{\psd}{PSD}        
\DeclareMathOperator{\nn}{NN}

\begin{document}

\title{Conic optimization for extremal geometry}

\address{F.~Vallentin, Department Mathematik/Informatik, Abteilung
  Mathematik, Universit\"at zu K\"oln, Weyertal~86--90, 50931 K\"oln,
  Germany}

\author{Frank Vallentin}
\email{frank.vallentin@uni-koeln.de}

\date{\today}

\subjclass{90C22, 52C17, 52C10, 46N10}

\date{October 1, 2025}

\maketitle

\begin{abstract} 
  The aim of this paper is to highlight recent progress in using conic
  optimization methods to study geometric packing problems. We will
  look at four geometric packing problems of different kinds: two on
  the unit sphere---the kissing number problem and measurable
  $\pi/2$-avoiding sets---and two in Euclidean space---the sphere
  packing problem and measurable one-avoiding sets.
\end{abstract}

\section{Introduction}

\subsection{Automatic reasoning for extremal problems in discrete geometry}

One central problem in discrete geometry is the \textit{kissing number
  problem}. The \textit{kissing number} in dimension $n$ is the
maximum number of non-overlapping unit spheres in $\R^n$ that can
simultaneously touch (``kiss'') a central unit sphere. The kissing
number in two dimensions obviously equals six, but the
three-dimensional case is already quite challenging. This case was
first studied by Isaac Newton and David Gregory in 1694 in connection
with the distribution of stars in the sky. In three dimensions, there
are infinitely many ways to arrange twelve unit spheres around a
central unit sphere. However, the question of whether a thirteenth
unit sphere can also touch the central one was only completely
resolved in 1953 by Sch\"utte and van der Waerden \cite{SchutteW1953}.

The thirteen-sphere problem can be expressed as a sentence in the
first-order theory of real closed fields:
\begin{equation}
\label{eq:thirteen-sphere}
\begin{split}
& \exists \;  x_{1,1}, x_{1,2}, x_{1,3}, \ldots, x_{13,1}, x_{13,2}, x_{13,3} : \\ 
& \qquad x_{i,1}^2 + x_{i,2}^2 + x_{i,3}^2 = 1 \text{ for } i = 1, \ldots, 13 \; \wedge \;  \\
& \qquad x_{i,1} x_{j,1} + x_{i,2} x_{j,2} + x_{i,3} x_{j,3} \leq 1/2 \text{ for } 1 \leq i < j \leq 13.
\end{split}
\end{equation}
Here the first set of equations ensures that each point
$(x_{i,1},x_{i,2},x_{i,3})$ lies on the unit sphere and represents the
contact point of the $i$-th sphere with the central unit sphere. The
second set of inequalities encodes that the angle between any two such
vectors is at least $\pi / 3$, which is equivalent to demanding that
the corresponding unit spheres do not overlap in their interiors.

It is a famous theorem of Tarski \cite{Tarski1948} (obtained around
1930 and published in 1948) that the first-order theory of real closed
fields is decidable, so there is an algorithm that decides whether
such a formula is true or false. However, Tarski's algorithm is not
practical, and more efficient algorithms were developed later.  (We
refer to the book of Basu, Pollack, and Roy \cite{BasuPR2006} for an
introduction to algorithmic real algebraic geometry.) For instance, it
is known that one can solve the existential theory of the reals
(sentences where all variables are bound to an $\exists$-quantifier,
like in the thirteen-sphere problem) using only polynomial space, thus
in exponential time.

Applying any of these methods to the thirteen-sphere problem is still
far beyond reach for current computers. (In fact the statement~\eqref{eq:thirteen-sphere} is
false: The kissing number in dimension three is twelve. See
Table~\ref{table:kissing-numbers} for the best upper and lower bounds
known for the kissing number in dimensions up to 24. Sch\"utte and van
der Waerden \cite{SchutteW1953} gave a classical proof applying
combinatorial and geometric arguments to show that there is no
thirteenth sphere.)

\smallskip

Nevertheless, the idea of using automatic reasoning to tackle extremal
problems in discrete geometry is both highly appealing and, by now,
widely established. Perhaps the most famous example is Hales’
resolution~\cite{Hales2005} of the Kepler conjecture, which asks for
the highest possible density of sphere packings in three
dimensions. His solution combined deep mathematical insight with an
elaborate computer-assisted argument, which was later fully formalized
and verified~\cite{Hales2017}.  Another, slightly different, line of
work uses optimization techniques to systematically search for
non-constructive bounds. A pioneering contribution in this direction
was the invention of the linear programming method of
Delsarte~\cite{Delsarte1973}. The power of this approach was
demonstrated spectacularly in Viazovska’s breakthrough
solution~\cite{Viazovska2017} of the sphere packing problem in
dimension 8.

\smallskip

The aim of this paper is to highlight recent progress in using conic
optimization methods to study geometric packing problems. Like
Tarski's algorithm, these methods have the potential to eventually
completely resolve such problems, provided sufficient computational resources are
available. There is, however, one essential difference from the
logical approach: optimization methods typically appear as a hierarchy
of increasingly tight relaxations. The first steps of this hierarchy
can be computed in practice, and each step may already lead to new and
interesting results. The main ingredients are conic optimization (in
particular semidefinite programming), symmetry reduction via harmonic
analysis, and techniques for rounding numerical approximations to
exact and easily verifiable solutions.

\subsection{Some extremal problems in discrete geometry}

Before turning to optimization methods, we briefly discuss the types
of extremal problems in discrete geometry to which they apply.

\smallskip

Many problems in discrete geometry are concerned with the optimal
distribution of finitely many points $\{x_1, \ldots, x_N\}$ in a
compact metric space $V$ equipped with a metric
$d$. There are many possibilities to define the quality of such a
configuration: One can maximize the packing density (or equivalently
the packing radius), which is by far the best-studied example. Other
important optimization problems include minimizing potential energy,
minimizing covering density, or the max-min polarization problem.

\begin{enumerate}
\item[(i)] \textit{Maximizing packing radius.}  How can we distribute
  $N$ points on the metric space $V$ so that the minimal
  distance between pairs of distinct points is maximized? In other
  words, we consider the optimization problem
\[
\max_{x_1, \ldots, x_N \in V} \;  \min\{d(x_i,x_j) : 1 \leq i < j \leq N\}.
\] 
This question, for example, is relevant to coding theory: One seeks to
distribute $N$ codewords in a manifold of possible signals
$V$ so as to minimize the probability of interference.

\item[(ii)] \textit{Minimizing potential energy.}  Given a potential
  function $p$, where $p(d(x,y))$ denotes the potential energy of two
  interacting particles $x, y \in V$, we consider the
  minimization problem
\[
\min_{x_1, \ldots, x_N \in V} \; \sum_{1 \leq i < j \leq N} p(d(x_i,x_j)).
\]
Such potential energy minimization problems arise naturally in the
study of physical particle systems. A classical example is the
\emph{Thomson problem}, which asks for the minimal-energy
configuration of $N$ points on the unit sphere $V = S^2$
interacting via the Coulomb potential $p(r) = 1/r$, where $r$ denotes
the Euclidean distance between two points.

Remarkably, certain highly symmetric configurations of a small number
of points are optimal for a broad class of natural potential
functions. For instance, the configuration of twelve points on $S^2$
forming the vertices of a regular icosahedron is optimal for many such
functions. This phenomenon, termed \emph{universal optimality}, was
identified by Cohn and Kumar~\cite{CohnK2007}.

Finally, maximizing the packing radius can be interpreted as a
limiting case of potential energy minimization when the potential
function is strictly decreasing in the distance and diverges as the
distance tends to zero.

\item[(iii)] \textit{Minimizing covering radius.}  How can we
  distribute $N$ points on the metric space $V$ so that the
  maximal distance to any other point on the metric space is minimized? In
  other words, we consider the optimization problem
\[
\min_{x_1, \ldots, x_N \in V} \; \max_{y \in V} \; \min \; \{d(x_i, y) : i = 1, \ldots, N\} .
\]
The problem of minimizing the covering radius is fundamental in metric
geometry. Example applications of covering codes, like data compression or football pools
are explained in \cite[Chapter 1.2]{CohenHLL1997}.

\item[(iv)] \textit{Max-min polarization.}  Let $p$ be a potential
  function. We consider the inhomogeneous variant of minimizing
  potential energy, given by the optimization problem
\[
\max_{x_1, \ldots, x_N \in V} \; \min_{y \in V} \; \sum_{i=1}^N p(d(x_i,y)).
\]
A physical interpretation of the inner minimization problem, proposed
by Borodachov, Hardin, and Saff \cite[Chapter 14]{BorodachovHS2019},
is as follows: If $p(d(x,y))$ represents the amount of a substance
received at a point $y$ due to an injector located at $x$, which
points receive the least substance when injectors are placed at $x_1,
\ldots, x_N$?

Analogous to the relationship between potential energy minimization
and packing radius maximization, max–min polarization can be viewed as
a limiting case of covering radius minimization.
\end{enumerate}

These geometric optimization problems have the flavor of binary $0/1$
optimization problems, which occur frequently in classical
combinatorial optimization: For every point $x \in V$ one has to make
the binary decision whether $x$ is chosen or not.

On the one hand, the geometric setting is more difficult than the
classical combinatorial setting, since the compact metric space $V$
may contain infinitely many points. Thus, one has to work with infinitely
many binary decision variables and the optimization problems become
infinite-dimensional. On the other hand, the geometric setting also
has advantages: Usually the geometric structure of $V$ is nice---it
is smooth and it has many symmetries---and one can exploit this when
performing the numerical optimization.

\subsection{Structure of the remainder of the paper}

In the following, we explain how tools from finite-dimensional
combinatorial optimization, particularly conic optimization
approaches, can be generalized to this infinite-dimensional geometric
setting.

In this paper, we focus on geometric packing problems, like maximizing
the packing radius. Related techniques for energy minimization and for
covering problems have also been investigated. We refer
to~\cite{Laat2020} for energy minimization and to~\cite{RienerRV2025}
for covering problems.

In Section~\ref{sec:Modeling}, we discuss how to model geometric
packing problems as independence numbers of graphs. Finding the
independence number of a given graph is a standard, though difficult
NP-hard problem in combinatorial optimization. We will look at four
geometric packing problems of different kinds: two on the unit
sphere---the kissing number problem and measurable $\pi/2$-avoiding
sets---and two in Euclidean space---the sphere packing problem and
measurable one-avoiding sets.

For highly structured graphs, in particular those with significant
symmetry, conic optimization approaches for determining the
independence number are known to perform best.  The \textit{conic
  optimization problems} we will mainly be concerned with are
optimization problems over convex cones of symmetric matrices; one
maximizes or minimizes a linear function over a convex cone
$\mathcal{K}$ intersected with an affine subspace. More precisely, we
consider the primal conic optimization problem
\begin{equation}
\label{eq:primal}
\begin{split}
p^* \; = \; \text{maximize} \quad & \langle C, X \rangle\\
\text{such that} \quad & X \in \mathcal{K},\\
& \langle A_j, X \rangle = b_j \; (j = 1, \ldots, m),
\end{split}
\end{equation}
where $C$, $A_1, \ldots, A_m$ are given symmetric matrices and $b_1,
\ldots, b_m \in \R$ are given real numbers. By $\langle A_j, X \rangle
= \mathrm{tr} (A_j X)$ we denote the trace inner product of symmetric
matrices (sometimes also called the Frobenius inner product). The
constraint $X \in \mathcal{K}$ is crucial here: it says that the
optimization variable, the symmetric matrix $X$, lies in the cone
$\mathcal{K}$. In this paper, the cones of positive semidefinite
matrices and the cone of completely positive matrices are used for
$\mathcal{K}$.  In the first case, we speak about \textit{semidefinite
  programming}, which is a vast, matrix-valued, generalization of
linear programming.

Conic optimization problems are convex optimization problems, so they
display a strong duality theory. The dual of~\eqref{eq:primal} is the
minimization problem
\begin{equation}
\label{eq:dual}
\begin{split}
d^* \; = \; \text{minimize} \quad & \sum_{j=1}^m b_j y_j\\
\text{such that} \quad & y_1, \ldots, y_m \in \R,\\
& \sum_{j=1}^m y_j A_j - C \in \mathcal{K}^*,
\end{split}
\end{equation}
where $\mathcal{K}^* = \{Y : \langle X, Y \rangle \geq 0 \text{ for
  all } X \in \mathcal{K}\}$ is the \textit{dual cone} of
$\mathcal{K}$. The cone of positive semidefinite matrices is
self-dual, but the dual cone of completely positive matrices is not
self-dual; it is the cone of copositive matrices. Weak duality $p^*
\leq d^*$ always holds between~\eqref{eq:primal} and \eqref{eq:dual}
and we have strong duality $p^* = d^*$ under some extra assumptions,
like strict feasibility (the existence of feasible solutions which lie
in the interior of the cones). Under mild technical assumptions,
semidefinite programs can be solved in polynomial time, in the sense
that the optimum can be approximated to within any desired precision
using the ellipsoid method~\cite{GrotschelLS1988} or interior-point
methods~\cite{KlerkV2016}. In practice, there are many implementations
of interior-point algorithms available. We refer to \cite{BenTalN2001,
  Nemirovski2007, LaurentV2026} for more details about conic
optimization and especially about the theory of semidefinite
programming.

In Section~\ref{sec:Formulations}, we derive conic optimization
formulations of the independence numbers for the graphs introduced in
Section~\ref{sec:Modeling}. These formulations are
infinite-dimensional analogues of the classical conic optimization
formulations of finite graphs.

Section~\ref{sec:Computations} describes how to solve these new conic
optimization problems. This is a highly nontrivial computational task
because the formulations are infinite-dimensional. We explain how to
exploit symmetry to simplify the computations and how to round
numerical solutions to obtain exact solutions in order to rigorously certify
the results.

In Section~\ref{sec:Results}, we give a survey of the results obtained
by this methodology and also discuss directions for future research.

\section{Modeling geometric packing problems as independence numbers of graphs}
\label{sec:Modeling}

Let~$G = (V, E)$ be an undirected graph (without loops and parallel
edges). A set~$I \subseteq V$ is \textit{independent} if it does not
contain pairs of adjacent vertices, that is, if for all~$x$, $y \in I$
we have~$\{x, y\} \notin E$. The \textit{independence number} of~$G$,
denoted by~$\alpha(G)$, is the maximum cardinality of an independent
set in~$G$. Complementary to the independence number is the
\textit{clique number}, which is the maximum cardinality of a
\textit{clique}, i.e., a set of pairwise adjacent vertices.

To model geometric packing problems as the independence number of a
graph, we extend the concept of independence number from finite to
infinite graphs. In this setting, the nature of both the vertex and
edge sets plays an essential role. Note also that in this model,
we fix the packing radius and maximize the number of points that can be placed.

Let~$V$ be a metric space with metric~$d$ and take~$D \subseteq (0,
\infty)$. The \textit{$D$-distance graph} on~$V$ is the graph~$G(V,
D)$ whose vertex set is~$V$ and in which vertices $x$,~$y$ are
adjacent if~$d(x, y) \in D$. Independent sets in~$G(V, D)$ are
sometimes called \textit{$D$-avoiding sets}. Let us consider a few
concrete choices for~$V$ and~$D$, corresponding to central problems in
discrete geometry. By $S^{n-1} = \{x \in \R^n : x \cdot x = 1\}$ we
denote the unit sphere with the Euclidean inner product~$x \cdot
y$. On the unit sphere we use the metric $d(x, y) = \arccos x \cdot
y$, the angle between the vectors~$x$ and~$y$.

\begin{enumerate}
\item[(i)] \textit{The kissing number problem: $V = S^{n-1}$
and~$D = (0, \pi / 3)$.} 
  
In this case, all independent sets in~$G(V, D)$ are finite; indeed,
also the independence number is finite. The independent sets in~$G(V,
D)$ are exactly the contact points of kissing configurations
in~$\R^n$, so~$\alpha(G(V, D))$ is the kissing number of~$\R^n$.

\item[(ii)] \textit{Measurable $\pi / 2$-avoiding sets: $V = S^{n-1}$
and~$D = \{\pi/2\}$.} 

An independent set in~$G(V, D)$ is a set without pairs of orthogonal
vectors. These sets can be infinite and even have positive surface
measure.  The right concept in this case is the \textit{measurable
  independence number}
\[
\alpha_\omega(G(V, D)) = \sup\{\, \omega(I) : \text{$I \subseteq V$ is
  measurable and independent}\,\},
\]
where~$\omega$ is the (normalized) surface measure on the sphere.

\item[(iii)] \textit{The sphere-packing problem: $V = \R^n$
and~$D = (0, 1)$.} 

Here we consider the Euclidean metric. The independent sets in~$G(V,
D)$ are the sets of centers of spheres in a packing of spheres of
radius~$1/2$ in~$\R^n$. So independent sets in~$G(V, D)$ can be
infinite but are always discrete and have Lebesgue measure~$0$. The
right definition of independence number in this case is the
\textit{sphere packing density}, informally the fraction of space
covered by the balls in the packing.  More precisely, we define the
\textit{upper density} of a Lebesgue-measurable set~$X \subseteq \R^n$
by
\[
\bar{\delta}(X) = \sup_{p\in\R^n}\limsup_{T\to\infty} \frac{\vol(X \cap (p +
  [-T, T]^n))}{\vol [-T,T]^n},
\]
where~$\vol$ is the Lebesgue measure. Then the \textit{sphere packing
  density} is
\[
\alpha_{\Delta}(G(V, D)) = \sup \{\, \bar{\delta}(I + 1/2 B_n) :
\text{$I \subseteq \R^n$ is independent}\},
\]
where $B_n = \{x \in \R^n : x \cdot x \leq 1\}$ is the unit ball.

\item[(iv)] \textit{Measurable one-avoiding sets: $V = \R^n$
  and~$D = \{1\}$.}
  
In this case, $G(V, D)$ is called the unit-distance graph
of~$\R^n$. Independent sets in this graph can be infinite and even
have infinite Lebesgue measure. So the right notion of independence
number is the \textit{independence density}, informally the fraction
of space covered. The \textit{independence density} is
\[
\alpha_{\bar{\delta}}(G(V , D)) = \sup \{\, \bar{\delta}(I) : \text{$I \subseteq \R^n$ is
  Lebesgue-measurable and independent}\}.
\]
\end{enumerate}

In the first two examples above, the vertex set is compact. For the
kissing number problem, there exists~$\delta > 0$ such that~$(0,
\delta) \subseteq D$. Then every point has a neighborhood that is a
clique (i.e., a set of pairwise adjacent vertices), which implies that
all independent sets are discrete and hence finite, given the
compactness of~$V$. For the second example,~0 is isolated
from~$D$. Then every point has an independent neighborhood and there
are independent sets of positive measure.

In the last two examples, the vertex set is not compact. For the
sphere packing problem, again there is~$\delta > 0$ such that~$(0,
\delta) \subseteq D$, and this implies that all independent sets are
discrete; since~$V$ is not compact, they can be infinite. For the
fourth example,~0 is again isolated from~$D$, hence there are
independent sets of positive measure and even infinite measure, given
that~$V$ is not compact.

We therefore see two factors at play. First, compactness of the vertex
set. Second, the nature of the edge set, which in the examples above
depends on~0 being isolated from~$D$ or not.

The graphs in examples~(i) and~(iii) are \textit{topological packing
  graphs}, a concept introduced by de Laat and
Vallentin~\cite{LaatV2015}. These are graphs in which the vertex set
carries a topology such that every finite clique is a subset of an
open clique. In particular, every vertex has a neighborhood that is a
clique.

The graphs in examples~(ii) and~(iv) are \textit{locally independent
  graphs}, which may be seen as the complements of topological packing
graphs. A topological graph is \textit{locally independent} if every
compact independent set is a subset of an open independent set.  In
particular, every vertex of a locally independent graph has an
independent neighborhood. The concept of locally independent graphs
was introduced by DeCorte, Oliveira, and
Vallentin~\cite{DeCorteOV2022}.

\section{Formulations and hierarchies of relaxations for the independence number}
\label{sec:Formulations}

The problem of computing the independence number of a finite graph is
NP-hard; in fact, its complementary problem, the maximum-clique
problem, appears in Karp's original list of~21 NP-hard problems. So it
is of interest to find good upper bounds which can be computed
efficiently. Starting from the Lov\'asz theta number, which is a
semidefinite programming relaxation of the independence number, we
describe two formulations of the independence number using conic
optimization. These lead to systematic methods that produce a
hierarchy of increasingly tight relaxations, eventually determining
the independence number. Furthermore, they serve as inspiration for
defining hierarchies for topological packings graphs and locally
independent graphs.

\subsection{The Lov\'asz theta number}

One of the best polynomial-time-computable upper bounds for the
independence number of a finite graph is the theta number, a graph
parameter introduced by Lovász~\cite{Lovasz1979} to determine the
Shannon capacity $\Theta(C_5)$ of the $5$-cycle graph $C_5$. Let~$G =
(V, E)$ be a finite graph. The theta number and its variants can be
defined in terms of the following conic optimization problem, in which
a linear function is maximized over the intersection of a convex cone
with an affine subspace:
\begin{equation}
\label{eq:finite-conic-prog}
\begin{split}
\vartheta(G, \mathcal{K}(V)) \; = \; \text{maximize} \quad & \langle J, A\rangle\\
\text{such that} \quad & \text{tr } A = 1,\\
& A(x, y) = 0 \quad \text{if~$\{x, y\} \in E$},\\
& A \in \mathcal{K}(V).
\end{split}
\end{equation}
Here, $A\colon V \times V \to \R$ is the optimization variable,
$J\colon V \times V \to \R$ is the all-ones matrix, $\langle J,
A\rangle = \text{tr } JA = \sum_{x,y \in V} A(x, y)$,
and~$\mathcal{K}(V) \subseteq \R^{V \times V}$ is a convex cone of
symmetric matrices.

The \textit{theta number} of~$G$, denoted by~$\vartheta(G)$, is
simply~$\vartheta(G, \psd(V))$, where~$\psd(V)$ is the cone of
positive semidefinite matrices where rows and columns are indexed by
the vertex set~$V$. In this case the conic optimization problem
becomes a semidefinite program, whose optimal value can be computed in
polynomial time. We have moreover $\vartheta(G) \geq \alpha(G)$: if~$I
\subseteq V$ is a nonempty independent set and $\chi^I \colon V \to
\{0,1\}$ is its characteristic vector, then $A = |I|^{-1} \chi^I
\otimes \chi^I$, which is the matrix such that
\[
A(x, y) = |I|^{-1} \chi^I(x) \chi^I(y),
\]
is a feasible solution of~$\vartheta(G, \psd(V))$; moreover~$\langle
J, A\rangle = |I|$, and hence $\vartheta(G) \geq |I|$. Since~$I$ is
any nonempty independent set, $\vartheta(G) \geq \alpha(G)$
follows. The theta number is a relaxation of the independence number
and it might happen (in fact it usually happens) that $\vartheta(G) >
\alpha(G)$. A strengthening of the Lov\'asz theta number is the
parameter $\vartheta'(G)$ introduced independently by McEliece,
Rodemich, and Rumsey~\cite{McElieceRR1978} and
Schrij\-ver~\cite{Schrijver1979}, obtained by taking $\mathcal{K}(V) =
\psd(V) \cap \nn(V)$, where~$\nn(V)$ is the cone of matrices with
nonnegative entries.

\subsection{A completely positive formulation}

Another choice for~$\mathcal{K}(V)$ is the cone
\[
\Ccal(V) = \cone\{\, f \otimes f : \text{$f\colon V \to \R$ and $f
  \geq 0$}\,\} \subseteq \psd(V) \cap \nn(V)
\]
of \defi{completely positive matrices}. The proof above that
$\vartheta(G) \geq \alpha(G)$ works just as well when $\Kcal(V) =
\Ccal(V)$. De Klerk and Pasechnik~\cite{KlerkP2002} observed that a
theorem of Motzkin and Straus~\cite{MotzkinS1965} implies that
$\vartheta(G, \mathcal{C}(V))$ equals $\alpha(G)$; a streamlined proof
of this fact goes as follows. If~$A$ is a feasible solution
of~$\vartheta(G, \Ccal(V))$, then, after suitable normalization,
\begin{equation}
\label{eq:A-convex-comb}
A = \alpha_1 f_1 \otimes f_1 + \cdots + \alpha_n f_n \otimes f_n,
\end{equation}
where~$\alpha_i > 0$, $f_i \geq 0$, and~$\|f_i\| = 1$ for
all~$i$. Since~$\|f_i\| = 1$, we have $\text{tr } f_i \otimes f_i = 1$, and
then since $\text{tr } A = 1$ we must have
$\alpha_1 + \cdots + \alpha_n = 1$. It follows that for some~$i$ we
have~$\langle J, f_i \otimes f_i\rangle \geq \langle J, A\rangle$;
assume then that this is the case for~$i = 1$.

Next, observe that since~$A(x, y) = 0$ for all~$\{x, y\} \in E$ and
each~$f_i$ is nonnegative, we must have~$f_1(x) f_1(y) = 0$ for
all~$\{x, y\} \in E$. This implies that~$I$, the support of~$f_1$, is an
independent set. Denoting by~$(f, g) = \sum_{x\in V} f(x) g(x)$ the
Euclidean inner product in~$\R^V$, we then have
\[
\langle J, A\rangle \leq \langle J, f_1 \otimes f_1\rangle = (f_1,
\chi_I)^2 \leq \|f_1\|^2 \|\chi_I\|^2 = |I| \leq \alpha(G)
\]
and, since~$A$ is any feasible solution, we get $\vartheta(G,
\Ccal(V)) \leq \alpha(G)$. Hence, $\vartheta(G,
\Ccal(V)) = \alpha(G)$.

\begin{theorem}
  Let $G = (V, E)$ be a finite graph. Then,
\begin{equation}
  \label{eq:ineq-chain}
  \vartheta(G,\psd(V)) \geq \vartheta(G, \psd(V) \cap \nn(V)) \geq
  \vartheta(G, \Ccal(V)) = \alpha(G).
\end{equation}
\end{theorem}

This seems to be mainly a curiosity: since solving $\vartheta(G,\Ccal(V))$ 
amounts to computing the independence number, computationally we have not 
gained anything. This is not entirely true, however: we now have a source of 
constraints that can be used to obtain better bounds. One such source of 
inequalities comes from the \emph{Boolean-quadratic cone}
\[
\mathrm{BQC}(V) = \cone\{\, f \otimes f : f \colon V \to \{0,1\}\,\} \subseteq \Ccal(V),
\]
which is a polyhedral cone. Valid inequalities of the Boolean-quadratic cone 
have been extensively studied, and many results are known; we refer to the book 
by Deza and Laurent~\cite[Chapter 5]{DezaL1997}. We have 
\[
\vartheta(G, \mathrm{BQC}(V)) \leq \vartheta(G, \Ccal(V)) = \alpha(G),
\]
but also $\vartheta(G, \mathrm{BQC}(V)) \geq \alpha(G)$, and hence $\vartheta(G, \mathrm{BQC}(V)) = \alpha(G)$.

Therefore, valid inequalities of the Boolean-quadratic cone can be used to 
strength\-en the theta number $\vartheta(G, \psd(V))$, yielding an upper bound 
for $\alpha(G)$ that may be strictly stronger than $\vartheta(G,  \psd(V))$. This process 
can be iterated, with more and more constraints added to strengthen the bound.

\smallskip

DeCorte, Oliveira, and Vallentin \cite{DeCorteOV2022} generalized the
completely positive formulation of the independence number to compact
locally independent graphs.  Problem~\eqref{eq:finite-conic-prog} can
be naturally extended to infinite topological graphs, as we will see
now. Let $G = (V, E)$ be a topological graph where~$V$ is
compact,~$\omega$ be a Borel measure on~$V$, $J \in L^2(V \times V)$
be the constant~1 kernel, and~$\Kcal(V) \subseteq L_{\text{sym}}^2(V
\times V)$ be a convex cone of symmetric kernels. When~$V$ is finite
with the discrete topology and~$\omega$ is the counting measure, the
following optimization problem is
exactly~\eqref{eq:finite-conic-prog}:
\begin{equation}
\label{eq:theta-problem}
\begin{split}
\vartheta(G, \Kcal(V)) \; = \; \text{maximize} \quad &\langle J, A\rangle\\
\text{such that} \quad & \int_V A(x, x)\, d\omega(x) = 1,\\
& A(x, y) = 0 \quad \text{if~$\{x, y\} \in E$},\\
&\text{$A$ is continuous and~$A \in \mathcal{K}(V)$.}
\end{split}
\end{equation}

The problem above is a straightforward extension
of~\eqref{eq:finite-conic-prog}, except that instead of the trace of
the operator~$A$ we take the integral over the diagonal, and we require $A$ to be continuous.

As before, there are several convex cones that can be put in
place of~$\Kcal(V)$. One is the cone~$\psd(V)$ of \textit{positive
kernels}, where we say a symmetric kernel $A \in L_{\text{sym}}^2(V
\times V)$ is \textit{positive} if for
all~$f \in L^2(V)$ we have
\[
\int_V \int_V A(x, y) f(x) f(y)\, d\omega(x) d\omega(y) \geq 0.
\]
The next cone is the cone of \defi{completely positive kernels} on~$V$, namely
\begin{equation}
\label{eq:cp-def}
\Ccal(V) = \text{cl cone}\{\, f \otimes f : \text{$f \in L^2(V)$ and~$f
  \geq 0$}\,\},
\end{equation}
with the closure taken in the norm topology on~$L^2(V \times V)$, and
where~$f \geq 0$ means that~$f$ is nonnegative almost everywhere.
Note that $\Ccal(V) \subseteq \psd(V)$, and hence~$\vartheta(G,
\psd(V)) \geq \vartheta(G, \Ccal(V))$.

The last cone is the \textit{Boolean-quadratic cone}
\[
\begin{split}
\mathrm{BQC}(V) = \text{cl}\{\, A \in L^2(V \times V) \; : \; & \text{$A$ is continuous and}\\
& \text{$A[U] \in \mathrm{BQC}(U)$ for all finite~$U \subseteq V$}\,\},
\end{split}
\]
with the closure taken in the $L^2$-norm topology
and where by~$A[U]$ we denote the restriction of~$A$ to~$U \times U$.

Under some extra, technical assumptions on~$G$ and~$\omega$, one has
$\vartheta(G, \Ccal(V)) = \vartheta(G, \mathrm{BQC}(V)) = \alpha_\omega(G)$, as in the finite case.
The proof of this theorem (see \cite[Theorem 5.1 and Theorem 7.1]{DeCorteOV2022} for
the exact statement) is fundamentally the same as in the finite case;
here is an intuitive description.

For the inequality $\vartheta(G, \Ccal(V)) \geq \alpha_\omega(G)$, one
constructs a feasible solution of~\eqref{eq:theta-problem} from any
independent set $I$ of $G$. Here one has to approximate the
characteristic function $\chi_I$ of $I$ by a continuous function $f :
V \to [0,1]$ so that the kernel $A = \|f\|^{-2} f \otimes f$ is a
feasible solution of \eqref{eq:theta-problem} with objective value
$\langle J, A\rangle \geq \omega(I) - \epsilon$.

For the reverse inequality $\vartheta(G, \Ccal(V)) \leq
\alpha_\omega(G)$ there are two key steps in the proof for finite
graphs as given above.  First, the matrix~$A$ is a convex combination
of rank-one nonnegative matrices, as
in~\eqref{eq:A-convex-comb}. Second, this together with the
constraints of our problem implies that the support of each~$f_i$
in~\eqref{eq:A-convex-comb} is an independent set. Then the support of
one of the~$f_i$'s will give us a large independent set.

In the proof that~$\vartheta(G, \Ccal(V)) = \alpha_\omega(G)$ for an
infinite topological graph we will have to repeat the two steps
above. Now~$A$ will be a kernel, so it will not be in general a convex
combination of finitely many rank-one kernels as
in~\eqref{eq:A-convex-comb}; Choquet theory~(see e.g.\ \cite[Chapters
  8--11]{Simon2011}) will allow us to express~$A$ as a sort of convex
combination of infinitely many rank-one kernels. Next, it will not be
the case that the support of any function appearing in the
decomposition of~$A$ will be independent, but depending on some
properties of~$G$ and~$\omega$ we will be able to fix this by removing
from the support the measure-zero set consisting of all points that
are not density points.

As the distance graph $G(S^{n-1}, \{\pi/2\})$ is a compact locally
independent graph, which also satisfies the extra technical
assumptions, we get the following exactness result for the measurable
independence number:
\[
\begin{split}
\vartheta(G(S^{n-1}, \{\pi/2\}), \Ccal(S^{n-1})) & = \vartheta(G(S^{n-1}, \{\pi/2\}), \mathrm{BQC}(S^{n-1}))\\
& = \alpha_\omega(G(S^{n-1}, \{\pi/2\}))
\end{split}
\]

Castro-Silva \cite{CastroSilva2021} provided an alternative proof
of this identity. In his approach, a nearly optimal kernel $A$ is
approximated in the supremum norm by a finite-rank completely positive
kernel $\tilde{A}$. However, after passing to such a finite-rank
approximation, one can no longer guarantee that $\tilde{A}$ vanishes
on the edges of $G(S^{n-1}, {\pi/2})$.  To overcome the errors
introduced by this approximation, Castro-Silva employs a
\textit{supersaturation} argument---a concept in extremal graph
theory, here adapted to the measurable setting. The key idea is that
if the objective value $\langle J, \tilde{A} \rangle$ of the
approximating kernel were significantly larger than the measurable
independence number, then the average value of $\tilde{A}$ on the
edges would necessarily be bounded away from $0$. This, however, would
contradict the fact that $\tilde{A}$ closely approximates $A$ in the
supremum norm.

One can also determine the independence density
$\alpha_{\bar{\delta}}(G(\R^n, \{1\}))$ of the unit-distance graph
$G(\R^n, \{1\})$ using a completely positive formulation. However,
this requires to work with a different cone of completely positive
functions on $\R^n$, which takes into account the
translation-invariance of the graph $G(\R^n, \{1\})$. This was done by
DeCorte, Oliveira, and Vallentin in \cite{DeCorteOV2022}.

A function~$f \in L^\infty(\R^n)$ is said to be of \textit{positive
  type} if~$f(x) = \overline{f(-x)}$ for all~$x \in \R^n$ and if for
every~$\rho \in L^1(\R^n)$ we have
\[
\int_{\R^n} \int_{\R^n} f(x-y) \rho(x) \overline{\rho(y)}\, dx dy \geq
0.
\]
The set of all
functions of positive type is a closed and convex cone, which we
denote by~$\psd(\R^n)$. A continuous function of positive
type~$f\colon\R^n \to \C$ has a well-defined \textit{mean value}
\[
M(f) = \lim_{T \to \infty} \frac{1}{\vol [-T,T]^n} \int_{[-T,T]^n}
f(x)\, dx.
\]
We define the cone of \textit{completely positive functions} on~$\R^n$, namely
\[
\begin{split}
  \Ccal(\R^n) = \mathrm{cl}\{\, f \in L^\infty(\R^n) \; : \; & \text{$f$ is continuous and}\\
  & \text{ $\bigl(f(x-y)\bigr)_{x,y\in U} \in \Ccal(U)$ for all
    finite~$U \subseteq \R^n$}\,\},
\end{split}
\]
where the closure is taken in the~$L^\infty$ norm; note
that~$\Ccal(\R^n)$ is a cone contained in~$\psd(\R^n)$. Finally, we define the cone of \textit{Boolean-quadratic functions} on~$\R^n$ by 
\[
\begin{split}
  \mathrm{BQC}(\R^n) = \text{cl}\{\, f \in L^\infty(\R^n) \; : \; & \text{$f$ is real
    valued and continuous and}\\
  & \text{$\bigl(f(x- y)\bigr)_{x,y \in U} \in \mathrm{BQC}(U)$ for all
    finite~$U \subseteq \R^n$}\,\},
\end{split}
\]
with the closure taken in the~$L^\infty$ norm. Note that $\mathrm{BQC}(\R^n)
\subseteq \Ccal(\R^n)$.

Let~$D \subseteq (0, \infty)$ be a set of forbidden distances and
$\Kcal(\R^n) \subseteq \psd(\R^n)$ be a convex cone; consider the
optimization problem
\begin{equation}
\label{eq:rn-cp-problem}
\begin{split}
\vartheta(G(\R^n, D), \Kcal(\R^n)) \; = \;  \text{maximize} \quad &M(f)\\
\text{such that}  \quad & f(0) = 1,\\
  & f(x) = 0 \quad \text{if~$\|x\| \in D$},\\
  & \text{$f\colon\R^n \to \R$ is continuous and~$f
      \in \Kcal(\R^n)$.}
\end{split}
\end{equation}

The bound $\vartheta(G(\R^n, D), \psd(\R^n))$ was introduced by
Oliveira and Vallentin~\cite{OliveiraV2010}. For the other choices
$\Kcal(\R^n) = \Ccal(\R^n)$ or $\Kcal(\R^n) = \mathrm{BQC}(\R^n)$, DeCorte, Oliveira, and
Vallentin~\cite{DeCorteOV2022} proved the following exactness result.

\begin{theorem}[Theorem~6.3 and Theorem~7.3 in \cite{DeCorteOV2022}]
\label{thm:rn-exactness}  
If~$D \subseteq (0, \infty)$ is closed, then
\[
\vartheta(G(\R^n, D),
\Ccal(\R^n)) = \vartheta(G(\R^n, D),
\mathrm{BQC}(\R^n)) = \alpha_{\bar{\delta}}(G(\R^n, D)).
\]
\end{theorem}

In particular, for measurable one-avoiding sets, 
\[
\vartheta(G(\R^n, \{1\}),
\Ccal(\R^n)) = \vartheta(G(\R^n, \{1\}),
\mathrm{BQC}(\R^n)) = \alpha_{\bar{\delta}}(G(\R^n, \{1\})).
\]

\subsection{A semidefinite programming hierarchy}

Another way to systematically obtain stronger bounds is to use the
\textit{Lasserre hierarchy} for $0/1$ polynomial optimization
problems. This hierarchy consists of a sequence of semidefinite
programs of growing size, whose optimal values converge to the the
optimal value of the original $0/1$ polynomial optimization problem.
The Lasserre hierarchy was introduced by Lasserre in
\cite{Lasserre2001}. He proved that it converges in finitely many
steps using Putinar's Positivstellensatz \cite{Putinar1993}, a
powerful result in real algebraic geometry. Shortly thereafter,
Laurent \cite{Laurent2003} provided a combinatorial proof, which we
use as a blueprint.

The definition of the Lasserre hierarchy requires some notation.  Let
$G = (V, E)$ be a graph with $n$ vertices. Let $t$ be an integer with
$0 \leq t \leq n$.  By $\mathcal{I}_t$ we denote the set of all
independent sets of $G$ of cardinality at most $t$.  A vector $y
\in \mathbb{R}^{\mathcal{I}_{2t}}$ defines a \textit{combinatorial
  moment matrix} of order $t$ by
\[
M_t(y) \in \psd(\mathcal{I}_t) \quad \text{with} \quad
(M_t(y))(J,J') = 
\begin{cases}
y({J \cup J'}) & \text{if } J \cup J' \in \mathcal{I}_{2t},\\
0 & \text{otherwise.}
\end{cases}
\]

For example, for the graph with vertex set $V = \{1,2,3\}$ and edge
set $E = \{\{1,2\}, \{1,3\}\}$, the combinatorial moment matrices of
order one and two have the following form:
\[
M_1(y) =
\bordermatrix{
& \emptyset & 1 & 2 & 3 \cr
\emptyset & y_{\emptyset} & y_1 & y_2 & y_3 \cr
1 & y_1 & y_1 & 0 & 0 \cr
2 & y_2 & 0 & y_{2} & y_{23} \cr
3 & y_3 & 0 & y_{23} & y_{3}\cr
},
\quad
M_2(y) =
\bordermatrix{ 
  & \emptyset & 1 & 2 & 3 & 23 \cr
\emptyset & y_{\emptyset} & y_1 & y_2 & y_3 & y_{23} \cr
1 & y_1 & y_1 & 0 & 0 & 0\cr
2 & y_2 & 0 & y_{2} & y_{23} & y_{23}\cr
3 & y_3 & 0 & y_{23} & y_{3} & y_{23}\cr
23 & y_{23} & 0 & y_{23} & y_{23} & y_{23}\cr
}.
\]
Here and in the following, we simplify notation and use $y_i$ instead
of $y(\{i\})$ and $y_{12}$ instead of $y({\{1,2\}})$.  Note that
$M_1(y)$ occurs as a principal submatrix of $M_2(y)$.

Let $t$ be an integer
with $1 \leq t \leq n$. The  \textit{Lasserre bound of $G$ of step
 $t$} is the value of the semidefinite program
\begin{equation*}
\begin{split}
\mathrm{las}_t(G) \; = \; \text{maximize} \quad &  \sum\nolimits_{i \in V} y_i \\
\text{such that} \quad & y \in \mathbb{R}^{\mathcal{I}_{2t}}_+,\\
& y_{\emptyset} = 1,\\
&  M_t(y) \in \psd(\mathcal{I}_t).
\end{split}
\end{equation*}

The Lasserre bounds form a hierarchy of stronger and stronger upper
bounds for the independence number of $G$ starting with
$\mathrm{las}_1(G) = \vartheta'(G)$, the strengthening of the Lov\'asz
theta number. Inequality $\mathrm{las}_{t+1}(G) \leq
\mathrm{las}_t(G)$ holds, because if $M_{t+1}(y)$ is positive
semidefinite, then also $M_{t}(y)$, being a principal submatrix of
$M_{t+1}(y)$, is positive semidefinite. Inequality $\mathrm{las}_t(G)
\geq \alpha(G)$ holds, because for any independent set $I$ of $G$ the
characteristic vector $\chi^I_{2t} \in \R^{\mathcal{I}_{2t}}$ defined
by
\[
\chi^I_{2t}(J) = 
\begin{cases}
  1 & \text{if } J \subseteq I,\\
  0 & \text{otherwise},
\end{cases}
\]
is a feasible solution of $\mathrm{las}_t(G)$.

Often one is interested in certifying an upper bound for
$\alpha(G)$. This can be done by exhibiting any feasible solution of
the conic optimization dual of the Lasserre bound of step $t$. The
conic optimization dual is a minimization problem, namely
\begin{equation}
\label{eq:Lasserre-bound-dual}
\begin{split}
\mathrm{las}_t(G) \; = \; \text{minimize} \quad & A(\emptyset, \emptyset) \\
\text{such that} \quad & A \in \psd(\mathcal{I}_t),\\
& \sum_{J, J' \in \mathcal{I}_t, J \cup J' = S} A(J,J')\\
& \quad \leq 
\begin{cases}
-1 & \text{if } S = \{i\} \text{ for some } i \in V,\\
0 & \text{if } S \in \mathcal{I}_{2t} \setminus
(\{\emptyset\} \cup \{\{i\} : i \in V\}).
\end{cases}
\end{split}
\end{equation}
Indeed, weak duality implies that any feasible solution $A$ of the
dual satisfies $A(\emptyset, \emptyset) \geq \mathrm{las}_t(G) \geq
\alpha(G)$. One can also directly verify the inequality $A(\emptyset,
\emptyset) \geq \alpha(G)$, which is crucial in applications, as
follows: If $I$ is an independent set of $G$, then
\[
0 \leq \sum_{J, J' \in \mathcal{I}_t, J \cup J' \subseteq I} A(J, J')
= \sum_{S \in \mathcal{I}_{2t}, S \subseteq I} \sum_{J, J' \in \mathcal{I}_t, J \cup J' = S} A(J, J')
\leq A(\emptyset, \emptyset) - |I|.
\]

An important feature of the Lasserre bound is that it does not lose
information. If the step of the hierarchy is high enough, we can
exactly determine the independence number of $G$: For every graph $G$
the Lasserre bound of step $t = \alpha(G)$ is exact; that means
$\mathrm{las}_t(G) = \alpha(G)$ for every $t \geq \alpha(G)$. This is
a consequence of the M\"obius inversion formula for partially ordered
sets, see \cite{Laurent2003} or \cite{LaatV2015} for a proof.

However, in practice, only the first few steps of the Lasserre
hierarchy can be computed, since the size of the combinatorial moment
matrix is usually of order $\Theta(n^t)$. However, in some favorable
cases, already the first steps give excellent bounds. For example,
when $G$ is a perfect graph, then even $\vartheta(G) = \alpha(G)$
holds (\cite{GrotschelLS1988}).

De Laat and Vallentin \cite{LaatV2015} generalized the Lasserre bound
to compact topological packing graphs. Whereas on the primal,
maximization side, the Lasserre bound for topological packing graphs
is defined using measures as optimization variables, the dual,
minimization side, is very close to the finite case. The only
difference is that the cone of positive semidefinite matrices is
replaced by the cone of continuous positive kernels.  After this change,
one can define the Lasserre bound of step $t$ for a topological
packing graph $G = (V,E)$ exactly as
in~\eqref{eq:Lasserre-bound-dual}. One main result of \cite{LaatV2015}
is that the generalization still satisfies all the properties of the
finite case:

\begin{theorem}
\label{thm:lasserre-bound}  
If $G$ be a compact topological packing graph, then
\[
\mathrm{las}_1(G) \geq
\mathrm{las}_2(G) \geq \ldots \geq
\mathrm{las}_{\alpha(G)}(G) = \alpha(G).
\]
\end{theorem}

The distance graph $G(S^{n-1}, (0,\pi/3))$ is a compact topological
packing graph and so we obtain a hierarchy of increasingly tight upper
bounds for the kissing number problem. The first step coincides with
the \textit{Delsarte-Goethals-Seidel linear programming bound}
\cite{DelsarteGS1977} as first realized by Bachoc, Nebe, Oliveira, and
Vallentin \cite{BachocNOV2009}.

Cohn and Salmon~\cite{CohnS2021} defined the Euclidean limit of the Lasserre bound for the sphere-packing graph $G(\R^n, (0,1))$. They showed that the first step coincides with the \textit{Cohn–Elkies linear programming bound}~\cite{CohnE2003}, and that for each $t$, the Euclidean limit of the $t$-th step provides an upper bound on the sphere-packing density. Moreover, the bounds converge to the sphere-packing density as $t \to \infty$.

\section{Computations: Symmetry reduction and rigorous verifications}
\label{sec:Computations}

In this section, we explain how to explicitly compute the Lasserre
bound, or variants thereof, in the case of the kissing number
problem. Similar techniques are available for other packing
problems. We refer to \cite{CohnE2003, LaatOV2014, CohnLS2022} for
computing bounds for the sphere packing problem; to
\cite{BachocNOV2009, DeCorteP2016,DeCorteOV2022, BekkerKOV2023} for
computing bounds for measurable $\pi/2$-avoiding sets; and to
\cite{OliveiraV2010, DeCorteOV2022} for computing bounds for
measurable $1$-avoiding sets.

\subsection{Exploiting symmetry via harmonic analysis}

When a graph has infinitely many vertices, then computing any step in
the semidefinite optimization hierarchies is an infinite-dimensional
semidefinite program. In most cases, we do not know how to solve these
optimization problems by analytic means. So one has to use a computer
to determine an, at least approximate, optimal solution. Therefore a
systematic approach to approximate the infinite-dimensional
optimization problem by a sequence of finite-dimensional ones is
needed.

One approach would be to discretize the graph and use the
``classical'' hierarchies for finite graphs. However, this is usually
not a good idea, since by discretizing the graph one destroys the
symmetry of the situation.

Another approach, the one which we advocate here, is to first
transform the semidefinite program at hand to its Fourier domain
(e.g.\ we work with the space of Fourier coefficients) and then
perform the discretization in the Fourier domain.
Since in the Fourier domain the symmetries are particularly visible, the full symmetry of
the situation can be exploited.  For this we compute \textit{explicit}
parametrizations of invariant convex cone of positive 
kernels in terms of their Fourier coefficients.

It is a well-known fact that symmetries can be very beneficially
exploited when solving convex optimization problems. We refer to
Bachoc, Gijswijt, Schrijver, Vallentin~\cite{BachocGSV2012} for a
survey on how to treat invariant semidefinite programs. For example
in \eqref{eq:Lasserre-bound-dual} the orthogonal group $\mathrm{O}(n)$
naturally acts on the optimization problem giving
$\mathrm{las}_t(G(S^{n-1},(0,\pi/3)))$. So it suffices to restrict the
optimization variable $A$ to the convex cone of
$\mathrm{O}(n)$-invariant positive kernels, denoted by
\begin{equation*}
\begin{split}
\psd(\mathcal{I}_t)^{\mathrm{O}(n)} = \{A \in L_{\text{sym}}^2(\mathcal{I}_t \times
\mathcal{I}_t) \; : \; & \text{$A$ positive kernel, and}\\ 
& \text{$A(\gamma J,\gamma gJ') = A(J,J')$}\\
& \;\; \text{for all $\gamma \in
  \mathrm{O}(n)$ and all $J, J' \in \mathcal{I}_t$}\}.
\end{split}
\end{equation*}

Abstractly, let $\Gamma$ be a compact matrix group acting transitively on a set
$V$ so that $V$ carries a Haar measure $\mu$ so that $\mu(\gamma S) = \mu(S)$ for all $\gamma \in
\Gamma$ and all measurable sets $S \subseteq V$. The action of
$\Gamma$ on $V$ extends to an action on the $L^2$-space of complex-valued
continuous functions $L^2(V) = L^2(V, \mu)$ via $(\gamma,f)(x) \mapsto (\gamma f)(x) = f(\gamma^{-1}x)$.

We want to compute an explicit parametrization of the cone $\psd(V)^\Gamma$ of $\Gamma$-invariant positive kernels. For this we can use the following recipe which is based on the celebrated Peter-Weyl theorem connecting group representations with Fourier analysis.

To state the Peter-Weyl theorem we need some vocabulary: A subspace $S
\subseteq L^2(V)$ is called \textit{$\Gamma$-invariant} if
$\gamma S = S$ for all $\gamma \in \Gamma$, i.e.\ if for every $\gamma
\in \Gamma$ and for every $f \in S$ we have $\gamma f \in S$ as
well. A nonzero subspace $S$ is called \textit{$\Gamma$-irreducible}
if $\{0\}$ and $S$ are the only $\Gamma$-invariant subspaces of
$S$. Let $S$ and $S'$ be two invariant subspaces, a linear map $T : S
\to S'$ is called a \textit{$\Gamma$-map} if $T(\gamma f) = \gamma
T(f)$ for all $\gamma \in \Gamma$, and $f \in L^2(V)$. We say
that $S$ and $S'$ are \textit{$\Gamma$-equivalent} if there is a
bijective $\Gamma$-map between them. Now the Peter-Weyl theorem
together with Schur orthogonality states: All irreducible subspaces of
$L^2(V)$ are of finite dimension and the Hilbert space $L^2(V)$
decomposes orthogonally as a completed direct sum
\begin{equation*}
L^2(V) = \widehat{\bigoplus_{k = 0}^\infty} H_k, \quad \text{and} \quad H_k = \bigoplus_{i = 1}^{m_k} H_{k,i},
\end{equation*}
where $H_{k,i}$ is $\Gamma$-irreducible, and $H_{k,i}$ is $\Gamma$-equivalent to
$H_{k',i'}$ if and only if the first index coincides, i.e.\ $k = k'$. The dimension $h_k$ of $H_{k,i}$
is finite the multiplicity $m_k$ is bounded by $h_k$.
In other words, $L^2(V)$ has a complete orthonormal system $e_{k,i,l}$,
where $k = 0, 1, \ldots$, $i = 1, 2, \ldots, m_k$, $l = 1, \ldots,
h_k$ so that
\begin{enumerate}
\item the space $H_{k,i}$ spanned by $e_{k,i,1}, \ldots, e_{k,i,h_k}$
is $\Gamma$-irreducible,
\item the spaces $H_{k,i}$ and $H_{k',i'}$ are $\Gamma$-equivalent if and only
  if $k = k'$,
\item there are $\Gamma$-maps $\phi_{k,i} : H_{k,1} \to H_{k,i}$
  mapping $e_{k,1,l}$ to $e_{k,i,l}$.
\end{enumerate}

The complete orthonormal system $e_{k,i,l}$ of the Peter-Weyl theorem
is very useful to characterize $\Gamma$-invariant, positive
kernels. This is the content of the following theorem which
essentially is due to Bochner \cite{Bochner1941}.

\begin{theorem}
\label{thm:bochner}
  Let $e_{k,i,l}$ be a complete orthonormal system for
  $L^2(V)$ as above. Every $\Gamma$-invariant, positive kernel
  $A \in \psd(V)^{\Gamma}$ can be written as (with convergence in $L^2$)
\begin{equation}
\label{eq:bochnerrep}
  A(x,y) = \sum_{k = 0}^\infty \sum_{i,j = 1}^{m_k}
  f_{k,ij} \sum_{l = 1}^{h_k} e_{k,i,l}(x) \overline{e_{k,j,l}(y)}
   = \sum_{k = 0}^\infty \langle F_k, Z_k(x,y) \rangle,
\end{equation}
with $(F_k)_{ij} = f_{k,ij}$ and $(Z_k(x,y))_{ij} = \sum_{l =
  1}^{h_k} e_{k,i,l}(x) \overline{e_{k,j,l}(y)}$ and where every $F_k$
(a matrix-valued Fourier coefficients of $A$) is
Hermitian positive semidefinite.
\end{theorem}

A classical example of this characterization is due to
Schoenberg~\cite{Schoenberg1942} for the sphere $V = S^{n-1}$ and the orthogonal group $\Gamma = \mathrm{O}(n)$ acting naturally on $S^{n-1}$. Here all the pieces of the puzzle
fall most neatly into place as $m_k = 1$ for all $k$ and $H_k$ is the space
of homogeneous harmonic polynomials of degree $k$.  More precisely,
let $\mathrm{Pol}_{\leq d}(S^{n-1})$ be the space of real polynomial
functions on $S^{n-1}$ of degree at most $d$, then
\begin{equation}
\mathrm{Pol}_{\leq d}(S^{n-1}) = H_0^n \oplus H_1^n \oplus \cdots \oplus H_d^n,
\end{equation}
where $H_k^n$ is the $\mathrm{O}(n)$-irreducible space of homogeneous,
harmonic polynomials of degree $k$ in $n$ variables; the dimension of
these spaces is denoted by $h_k^n = \dim(H_k^n)$. Schoenberg's
characterization states that all $\mathrm{O}(n)$-invariant, continuous, positive type kernel on $S^{n-1}$ are of the form
\begin{equation}
\label{eq:Schoenberg}
\sum_{k=0}^{\infty} f_k P^n_k(x
\cdot y) \quad \text{with $f_k \geq 0$, $\sum_{k=0}^{\infty} f_k < \infty$},
\end{equation}
where $P_k^n$ is the polynomial of degree $k$ satisfying the
orthogonality relation
\begin{equation*}
\label{eq:orthogonality-relation-pnk}
\int_{-1}^{1} P_k^n(t) P_l^n(t) (1-t^2)^{\frac{n-3}{2}} dt = 0 \text{ if
}  k\neq l,
\end{equation*}
and where the polynomial $P_k^n$ is normalized by $P_k^n(1) = 1$.  The
polynomials $P_k^n$ appear under different names with different
normalizations: Jacobi polynomials, Gegenbauer polynomials,
ultraspherical polynomials are the most common ones.  The equality
in~\eqref{eq:Schoenberg} should be interpreted as follows: A kernel
$A \in L^2_{\text{sym}}(S^{n-1} \times S^{n-1})$ is $\mathrm{O}(n)$-invariant, continuous, and positive if and only if there are
nonnegative numbers $f_0, f_1, \ldots$ so that the series
$\sum_{k=0}^{\infty}f_k$ converges and so that
\[
A(x,y) = \sum_{k=0}^{\infty} f_k P^n_k(x \cdot y)
\]
holds. Here the right-hand side even converges absolutely and uniformly over $S^{n-1} \times S^{n-1}$.

Schoenberg's characterization is the basic technical tool to turn the
semidefinite program defining $\mathrm{las}_1(G(S^{n-1},(0,\pi/3)))$
into the Delsarte-Goethals-Seidel linear programming bound (note that $\mathcal{I}_1 = \{\{x\} : x \in S^{n-1}\}
\cup \{\emptyset\}$):
\begin{equation*}
\begin{split}
\mathrm{las}_1(G(S^{n-1},(0,\pi/3))) \; = \; \text{minimize} \quad & \lambda\\
\text{such that} \quad &  \lambda \in \R, \; f_0, f_1, \ldots \geq 0, \; \sum_{k = 0}^{\infty} f_k  = \lambda-1,\\
& \sum_{k = 0}^{\infty} f_k P_k^n(t)  \le  -1 \text{ for all } t \in [-1,1/2].
\end{split}
\label{eq:lp-bound}
\end{equation*}
This linear programming bound is also called a \textit{two-point
  bound} because it involves constraints on the two point distribution
of a configuration.

More complicated is the characterization of the cone
$\psd(S^{n-1})^{\mathrm{O}(n-1)}$. We consider $\mathrm{O}(n-1)$ as the subgroup of
$\mathrm{O}(n)$ which stabilizes one point $e \in S^{n-1}$ on the unit
sphere, the North pole. This falls slighty outside of the above mentioned recipe, since the action of $\mathrm{O}(n-1)$ on $S^{n-1}$ is not transitive. However, Bachoc and Vallentin~\cite{BachocV2009} showed that the recipe can be adapted to this situation. Under the action of $\mathrm{O}(n-1)$ we have the following decomposition:
\begin{equation*}
\mathrm{Pol}_{\leq d}(S^{n-1}) = \bigoplus_{k=0}^d \; (H_{k,k}^{n-1}
\oplus \cdots \oplus H_{k,d}^{n-1}),
\end{equation*}
where, for $i\geq k$, $H_{k,i}^{n-1}$ is the unique subspace of
$H_i^n$ isomorphic to $H_k^{n-1}$. Using the recipe one gets
\[
Z_k(x,y) = \big(Y_k^n\big)_{i,j}(e \cdot x, e \cdot y, x \cdot y),
\]
where we have, for all $0\leq i,j\leq d-k$,
\begin{equation}
\big(Y_k^n\big)_{i,j}(u,v,t) = \lambda_{i,j}
P_{i}^{n+2k}(u)P_{j}^{n+2k}(v)Q_k^{n-1}(u,v,t),
\end{equation}
and
\[
\displaystyle
Q_k^{n-1}(u,v,t) =\big((1-u^2)(1-v^2)\big)^{k/2}P_k^{n-1}\Big(\frac{t-uv}{\sqrt{(1-u^2)(1-v^2)}}\Big),
\]
and with normalization constants (recall that $e_{k,i,l}$ are orthonormal; in fact for our application the normalization is not crucial)
\[
\lambda_{i,j}=\frac{\omega_n}{\omega_{n-1}}\frac{\omega_{n+2k-1}}{\omega_{n+2k}}(h_i^{n+2k}h_j^{n+2k})^{1/2},
\]
where $\omega_n$ is the (standard Lebesgue, non-normalized) surface
area of $S^{n-1}$. 

The characterization of the cone $\psd(S^{n-1})^{\mathrm{O}(n-1)}$
leads to the \textit{three-point bound} for the kissing number by
Bachoc and Vallentin \cite{BachocV2008} where constraints on the three
point distribution are taken into account.  Set $S_k^n = \sum_{\sigma}
\sigma Y_k^n$, where $\sigma$ runs through all permutation of the
variables $u, v, t$. A simplified version of the three-point bound,
from \cite[Theorem 6.10]{BachocGSV2012}, is as follows:

\begin{equation*}
\begin{split}
\alpha(G(S^{n-1},(0,\pi/3))) \; \leq \; \text{minimize} \quad & 1 + \langle F_0, J_{d+1}\rangle\\
\text{such that} \quad & F_0 \in \psd(d+1), F_1 \in \psd(d), \ldots, F_d \in \psd(1) \\
& \displaystyle \sum_{k=0}^d \langle F_k, S_k^n (u,u,1)\rangle \leq
-\frac{1}{3}, \quad -1\leq u\leq 1/2\\
& \displaystyle \sum_{k=0}^d \langle F_k, S_k^n (u,v,t)\rangle \leq
0,\\
& \quad \text{$-1\leq u,v,t \leq 1/2,\; 1+2uvt-u^2-v^2-t^2\geq 0$,}
\end{split}
\end{equation*}
where $\psd(d+1), \psd(d), \ldots, \psd(1)$ denote the cones of
positive semidefinite matrices of sizes $(d+1) \times (d+1)$, $d
\times d$, \dots $1 \times 1$. This three-point bound is inspired by
Schrijver's \cite{Schrijver2005} three-point bound for binary error
correcting codes. Laurent~\cite{Laurent2007} showed that Schrijver's
three-point bound lies between the first and second step of the
Lasserre bound; the second step being a \textit{four-point bound}.

Very recently, de Laat, Leijenhorst, and de Muinck Keizer
\cite{LaatLM2024}, building on \cite{LaatMM2022}, were able to give an explicit, though heavily
computer assisted, parametrization of the cone
$\psd(\mathcal{I}_2)^{\mathrm{O}(n)}$, which made it possible to
compute the second step of the Lasserre bound
$\mathrm{las}_2(G(S^{n-1}, (0,\pi/3)))$.

\subsection{Computation and verification of bounds}

To compute bounds with the assistance of a computer, one must solve a
semidefinite program, which in primal standard form is given
in~\eqref{eq:primal} with the cone $\mathcal{K} = \psd(n)$ of positive
semidefinite matrices of size $n \times n$. To obtain good or even
optimal bounds, it is necessary to determine the optimal value $p^*$
with high accuracy, or even exactly. For this, one uses an
implementation of a semidefinite programming solver. One problem is
that all existing implementations that produce high accuracy solutions
are numerical interior-point solvers. These solvers ideally produce a
numerical approximation $X^*$ of a relative interior point of the
optimal face
\[
\mathcal{F} = \{Y \in \psd(n) : \langle C, Y \rangle = p^*, \langle A_j, Y \rangle = b_j \; (j = 1, \ldots, m)\}.
\]
(Under mild technical assumptions, interior-point algorithms, as they
follow the central path, converge to the analytic center of an optimal
face \cite{Halicka2002}.)  However, this means that $X^*$ is usually
neither positive semidefinite (it can have slightly negative
eigenvalues) nor does it satisfy the linear constraints $\langle A_j,
X \rangle = b_j$.

To address this issue, several rounding methods have been proposed and
implemented \cite{Monniaux2011, DostertLM2021, CohnLL2024}. The very first step of a
rounding method is to obtain a numerical approximation $X^*$, with
extremely high precision, of a relative interior-point of the optimal
face. For this, Leijenhorst and de Laat \cite{LeijenhorstL2024}
developed a high-precision primal-dual interior-point solver that
exploits additional low-rank structure to speed up the computation. In
the second step the numerical approximation $X^*$ is used to identify
the affine hull of the optimal face. It is known \cite{HillW1987} that
the minimal face of the cone $\psd(n)$ containing a given matrix $X
\in \psd(n)$ is
\[
\{Y \in \psd(n) : \ker X \subseteq \ker Y\}.
\]
Hence, if $X$ lies in the relative interior of the optimal face, then
all points $Y \in \mathcal{F}$ satisfy $\ker X \subseteq \ker Y$. To
detect the kernel, the LLL lattice basis reduction algorithm
\cite{LenstraLL} is used. In the third step, one performs a
\textit{facial reduction}, a coordinate transformation that transforms
the optimal face to become full-dimensional in a cone of positive
semidefinite matrices of smaller dimension. This transforms $X^*$ to
$\hat{X}^*$. In the last step $\hat{X}^*$ is rounded to the
transformed optimal face. After replacing each entry of $\hat{X}^*$ by
a close approximation in some fixed algebraic number field (usually
$\mathbb{Q}$ or $\mathbb{Q}(\sqrt{2})$ suffices), one obtains a point
in the transformed optimal face by exactly solving a least-squares
system.

It turns out that this rounding heuristic is quite successful and that
the numbers involved remain well-behaved. This contrasts with the
study of Nie, Ranestad, and Sturmfels~\cite{NieRS2010}, which
considers generic semidefinite programs with rational input. It would
be desirable to gain a deeper understanding of when and why the
rounding heuristic succeeds.

\section{Results and conclusion}
\label{sec:Results}

\subsection{The kissing number problem}

One highly influential and by now classical resource on the kissing
number problem is the book by Conway and Sloane~\cite{ConwayS1999}. In
the early years after its first edition in 1988, progress on improving
either lower or upper bounds was slow. It was widely believed that the
lower bounds reported there were in fact the correct values, and that
the available techniques for proving upper bounds—most notably the
Delsarte–Goethals–Seidel bound, used by Odlyzko and
Sloane~\cite{OdlyzkoS1979} and by Levenshtein~\cite{Levenshtein1979}
to solve the kissing number problem in dimensions 8 and 24—were not
strong enough to go further.

Over the past 20 years, beginning with Musin’s solution of the kissing
number problem in dimension 4~\cite{Musin2008}, first announced in
2004, the gap between lower and upper bounds for the kissing number
has steadily narrowed. This progress is due in large part to the
development of semidefinite programming bounds. At the same time, new
and sometimes surprising geometric constructions of spherical codes
have led to improved lower bounds. Remarkably, a recent improvement in
dimension 11 was achieved with the aid of artificial intelligence,
combining a large language model with genetic
programming~\cite{Novikov2025}.

In Table~\ref{table:kissing-numbers}, we provide an update of
\cite[Table~1.5]{ConwayS1999}, including references for the new
entries. For the most up-to-date records, see Cohn’s table of kissing numbers\footnote{\url{https://dspace.mit.edu/handle/1721.1/153312}}.

\begin{table}[htb]
\caption{Known bounds for the kissing number in various dimensions.
  This table updates \cite[Table~1.5]{ConwayS1999} of Conway and
  Sloane, with references provided for the new entries.}
\label{table:kissing-numbers}
\begin{center}
\begin{tabular}{cccc}
\toprule
Dimension $n$ & Lower bound & Upper bound & References\\
\midrule
3	& 12 & 12 &	\\
4	& 24 & 24 & \cite{Musin2008}, \cite{LaatLM2024}\\
5	& 40 & 44	& \cite{MittelmannV2010} \\
6	& 72 & 77 & \cite{LaatLM2024}\\
7	& 126 & 134	& \cite{MittelmannV2010}\\
8	& 240	& 240 & \\
9	& 306	& 363	&	\cite{MachadoO2018}\\
10 & 510 & 553 & \cite{Ganzhinov2025}, \cite{MachadoO2018} \\  
11 & 593 & 868 & \cite{Novikov2025} \cite{LeijenhorstL2024}\\
12 & 840 & 1355 & \cite{LeijenhorstL2024}\\
13 & 1154	& 2064 & \cite{Zinoviev1999}, \cite{LeijenhorstL2024}\\
14 & 1932	& 3174 & \cite{Ganzhinov2025} \cite{LeijenhorstL2024}\\
15 & 2564	& 4853 & \cite{LeijenhorstL2024}\\
16 & 4320	& 7320 & \cite{LeijenhorstL2024}\\
17 & 5730 & 10978	& \cite{CohnL2024}, \cite{LeijenhorstL2024}\\
18 & 7654 & 16406	& \cite{CohnL2024}, \cite{LeijenhorstL2024}\\
19 & 11692 & 24417 & \cite{CohnL2024}, \cite{LeijenhorstL2024}\\
20 & 19448 & 36195 & \cite{CohnL2024}, \cite{LeijenhorstL2024}\\
21 & 29768 & 53524 & \cite{CohnL2024}, \cite{LeijenhorstL2024}\\
22 & 49896 & 80810 & \cite{LeijenhorstL2024}\\
23 & 93150 & 122351 & \cite{LeijenhorstL2024}\\
24 & 196560	& 196560 & \\
\bottomrule
\end{tabular}
\end{center}
\end{table}

In dimension 4, Musin~\cite{Musin2008} showed that the kissing number
is $24$, using a combination of the Delsarte–Goethals–Seidel linear
programming bound with additional geometric arguments. More recently,
de Laat, Leijenhorst, and de Muinck Keizer~\cite{LaatLM2024} proved
that $\text{las}_2(G(S^3,(0,\pi/3))) = 24$. This allowed them to
establish that the configuration of $24$ points arising from the $D_4$
root system (or equivalently from the 24-cell) is unique, up to
orthogonal transformations. Most of the improvements in the upper bounds stem from the three-point
bound~\cite{BachocV2008}. Improved implementations of this bound led
to further progress~\cite{MittelmannV2010, MachadoO2018,
  LeijenhorstL2024}. One notable exception is dimension 6, where the
second step of the Lasserre hierarchy was applied in~\cite{LaatLM2024}
to surpass the three-point bound.

\smallskip

Looking ahead, further improvements of the upper bounds appear to
depend on implementing semidefinite programming bounds for higher
steps of the Lasserre hierarchy. This is highly nontrivial: performing
the symmetry reduction is already computationally demanding, and
solving the resulting semidefinite programs becomes increasingly
difficult. The rounding procedures required for rigorous bounds also
grow more involved. Nevertheless, one can envision the development of
a fully formal proof system capable of automatically verifying both
the numerical computations and the rounding steps.

\subsection{The sphere packing problem}

The sphere packing problem has a long history; we refer to
\cite{ConwayS1999} for further information. Since the publication of
\cite{ConwayS1999}, essentially all upper bounds for sphere packings in \cite[Table 1.2]{ConwayS1999}
have been improved. The Cohn–Elkies linear programming bound (being a two-point bound) played a
pivotal role in these improvements and, in particular, led to the
solution of the sphere packing problem in dimensions 8 and 24.

After the breakthrough results on the sphere packing problem
\cite{Viazovska2017, CohnKMRV2017}, the power of the Cohn–Elkies bound
is fairly well understood in low dimensions: it is known to give tight
bounds in dimensions 1, 8, and 24, and conjecturally also in dimension
2. In all other dimensions, the Cohn–Elkies bound is conjectured not
to be tight. Recently, Cohn, de Laat, and Salmon~\cite{CohnLS2022}
computed three-point bounds for $\alpha_\Delta(G(\R^n, (0,1)))$. This
provided new upper bounds for the sphere packing density in dimensions
4 through 7 and 9 through 16. For the most up-to-date records, see
Cohn’s table of sphere packing density bounds\footnote{\url{https://hdl.handle.net/1721.1/153311}}.

\smallskip

It is natural to ask how semidefinite programming bounds could be
used to prove tight results. In principle, they have this potential,
since they are known to converge. At present, however, there is no
numerical evidence indicating which steps would be required to establish
tightness in any dimension. Once such evidence becomes available, one
could attempt to adapt Viazovska’s techniques to the semidefinite
programming setting; for now, though, this seems out of reach.

\subsection{Measurable $\pi/2$-avoiding sets and the double cap conjecture}

The problem of determining the maximum surface measure of a
$\pi/2$-avoiding set was first posed by Witsenhausen~\cite{Witsenhausen1974}.
He obtained an upper bound of $1/n$ times the surface measure of the sphere
$S^{n-1}$ using a simple averaging argument, which is sharp for $S^1$.
Indeed, two antipodal open spherical caps of radius $\pi/4$ form a subset
with no pairs of orthogonal vectors. Kalai~\cite[Conjecture~2.8]{Kalai2015}
conjectured that this construction is optimal. This conjecture is known as
\emph{Kalai's double cap conjecture}, and it remains open for all $n \geq 3$.

The best upper bounds are due to Bekker, Kuryatnikova, Oliveira, and
Vera \cite{BekkerKOV2023}; see
Table~\ref{table:measurable-half-pi-avoiding}. Their computation is
based on the completely positive formulation of
$\alpha_\omega(G(S^{n-1},{\pi/2}))$. They develop a hierarchy of
semidefinite programs that approximate the completely positive cone,
which they further strengthen by using inequalities from the Boolean-quadratic cone.

\begin{table}[htb]
\caption{Bounds for the measurable independence number $\alpha_\omega(G(S^{n-1}, \{\pi/2\})$, where~$\omega$ is the (normalized) surface measure on the sphere. The lower bounds give the measure of a double cap. The upper bounds for $n \geq 3$ are all from \cite{BekkerKOV2023}.}
\label{table:measurable-half-pi-avoiding}
\begin{center}
\begin{tabular}{ccc}
\toprule
Dimension $n$ & Lower bound & Upper bound \\
\midrule
2 & 0.5 & 0.5 \\
3 & 0.2928\dots & 0.297742\\
4 & 0.1816\dots & 0.194297\\
5 & 0.1161\dots & 0.134588\\
6 & 0.0755\dots & 0.098095\\
7 & 0.0498\dots & 0.075751\\
8 & 0.0331\dots & 0.061178\\
\bottomrule
\end{tabular}
\end{center}
\end{table}

\subsection{Measurable one-avoiding sets and a conjecture by Erd\H{o}s}

The problem of finding measurable one-avoiding sets appears in Moser's
collection of problems~\cite{Moser1991}, partially compiled in 1966,
and was later popularized by Erd\H{o}s \cite{Erdoes1985}, who
conjectured that~$\alpha_{\bar{\delta}}(G(\R^2,\{1\})) < 1/4$
(cf.~Székely~\cite{Szekely2002}).

Another long-standing conjecture of Moser (cf.~Conjecture~1 in
Larman and Rogers~\cite{LarmanR1972}), related to Erd\H{o}s’s
conjecture, would imply that $\alpha_{\bar{\delta}}(G(\R^n,\{1\})) \leq
1/2^n$ for all $n \geq 2$. Moser’s conjecture asserts that the maximum
measure of a subset of the unit ball containing no pair of points at
distance~1 is at most $1/2^n$ times the measure of the unit ball. This
conjecture was shown to be false~\cite{OliveiraV2019}: the behavior of
subsets of the unit ball avoiding distance~1 resembles that predicted
by Kalai’s double cap conjecture.

To date, the best lower bound
$\alpha_{\bar{\delta}}(G(\R^2,\{1\})) \geq 0.22936$ is due to
Croft~\cite{Croft1967}, who placed \textit{tortoises} on the hexagonal
lattice. Here, a tortoise is defined as the intersection of an open
disc of radius $1/2$ with an open regular hexagon of height
$x = 0.96533\dots$. More recently, Ambrus, Csiszárik, Matolcsi, Varga,
and Zsámboki~\cite{AmbrusCMVZ2024} resolved Erd\H{o}s’s conjecture by
proving that $\alpha_{\bar{\delta}}(G(\R^2,\{1\})) < 0.2470$. Their
bound can be interpreted as arising from the completely positive
formulation of $\alpha_{\bar{\delta}}(G(\R^2,\{1\}))$, in which they
strengthened $\vartheta(G(\R^n, \{1\}), \psd(\R^n))$ with carefully chosen
inequalities from the Boolean-quadratic cone $\mathrm{BQC}(\R^n)$. The decisive step lay
in the selection of these inequalities; for this, they relied heavily
on massive computational power combined with a clever implementation
of a beam search algorithm.

\smallskip

For the last two problems, conic optimization has yielded the best
known upper bounds, but one might wonder whether this approach is
truly effective, as in no case---except for the trivial case of
$S^1$---is the bound tight. Identifying and analyzing tight cases in
this setting could provide valuable insight. In contrast,
comparatively little work has been done on establishing good lower
bounds for measurable one-avoiding sets.

\subsection{Beyond geometric graphs: Geometric hypergraphs}

In this paper, we focused on geometric packing problems that can be
formulated using the independence number of graphs. However, the
methods can be extended further to model packing problems via the
independence number of geometric hypergraphs. This provides a natural
framework for Euclidean Ramsey theory. The central question of
Euclidean Ramsey theory is: given a finite configuration $P$ of points
in $\R^n$ and an integer $r \geq 1$, does every $r$-coloring of $\R^n$
contain a monochromatic congruent copy of $P$? Conic optimization
methods for such questions have been studied in
\cite{CastroSilvaOSV2022, CastroSilvaOSV2023, CastroSilva2021,
  CastroSilva2023}, yielding improved bounds in Euclidean Ramsey
theory.

Let us end with a conjecture, stated in \cite[Conjecture~1]{CastroSilva2023} and also related to results of Bourgain~\cite{Bourgain1986} and Furstenberg, Katznelson, and Weiss~\cite{FurstenbergKW1990}, which falls within the framework of independence numbers of geometric hypergraphs but has so far resisted attack by conic optimization techniques: Let $A \subseteq \R^2$ be a set of positive upper density and let $u, v, w \in \R^2$ be noncollinear points. Then there exists $t_0 > 0$ such that for any $t \geq t_0$, the set $A$ contains a configuration congruent to $\{tu, tv, tw\}$.

\section*{Acknowledgments}

I am grateful to Davi Castro-Silva, Henry Cohn, David de Laat, Fernando Oliveira, and Marc Christian Zimmermann for their valuable feedback on the present paper.


\begin{thebibliography}{10}

\bibitem{AmbrusCMVZ2024}
{\sc G.~Ambrus, A.~Csiszárik, M.~Matolcsi, D.~Varga, and P.~Zsámboki}, {\em The density of planar sets avoiding unit distances}, Math. Program., 207 (2024), pp.~303--327, \url{https://doi.org/10.1007/s10107-023-02012-9}.

\bibitem{BachocGSV2012}
{\sc C.~Bachoc, D.~C. Gijswijt, A.~Schrijver, and F.~Vallentin}, {\em Invariant semidefinite programs}, in Handbook on semidefinite, conic and polynomial optimization, vol.~166 of Internat. Ser. Oper. Res. Management Sci., Springer, New York, 2012, pp.~219--269, \url{https://doi.org/10.1007/978-1-4614-0769-0\_9}.

\bibitem{BachocNOV2009}
{\sc C.~Bachoc, G.~Nebe, F.~M. de~Oliveira~Filho, and F.~Vallentin}, {\em Lower bounds for measurable chromatic numbers}, Geom. Funct. Anal., 19 (2009), pp.~645--661, \url{https://doi.org/10.1007/s00039-009-0013-7}.

\bibitem{BachocV2008}
{\sc C.~Bachoc and F.~Vallentin}, {\em New upper bounds for kissing numbers from semidefinite programming}, J. Amer. Math. Soc., 21 (2008), pp.~909--924, \url{https://doi.org/10.1090/S0894-0347-07-00589-9}.

\bibitem{BachocV2009}
{\sc C.~Bachoc and F.~Vallentin}, {\em Semidefinite programming, multivariate orthogonal polynomials, and codes in spherical caps}, European J. Combin., 30 (2009), pp.~625--637, \url{https://doi.org/10.1016/j.ejc.2008.07.017}.

\bibitem{BasuPR2006}
{\sc S.~Basu, R.~Pollack, and M.-F. Roy}, {\em Algorithms in real algebraic geometry}, vol.~10 of Algorithms and Computation in Mathematics, Springer-Verlag, Berlin, second~ed., 2006.

\bibitem{BekkerKOV2023}
{\sc B.~Bekker, O.~Kuryatnikova, F.~M. de~Oliveira~Filho, and J.~C. Vera}, {\em Optimization hierarchies for distance-avoiding sets in compact spaces}, 2023, \url{https://arxiv.org/abs/2304.05429}.

\bibitem{BenTalN2001}
{\sc A.~Ben-Tal and A.~Nemirovski}, {\em Lectures on modern convex optimization: Analysis, algorithms, and engineering applications}, MPS/SIAM Series on Optimization, Society for Industrial and Applied Mathematics (SIAM), Philadelphia, PA; Mathematical Programming Society (MPS), Philadelphia, PA, 2001, \url{https://doi.org/10.1137/1.9780898718829}.

\bibitem{Bochner1941}
{\sc S.~Bochner}, {\em Hilbert distances and positive definite functions}, Ann. of Math. (2), 42 (1941), pp.~647--656, \url{https://doi.org/10.2307/1969252}.

\bibitem{BorodachovHS2019}
{\sc S.~V. Borodachov, D.~P. Hardin, and E.~B. Saff}, {\em Discrete energy on rectifiable sets}, Springer Monographs in Mathematics, Springer, New York, 2019, \url{https://doi.org/10.1007/978-0-387-84808-2}.

\bibitem{Bourgain1986}
{\sc J.~Bourgain}, {\em A {S}zemer\'edi type theorem for sets of positive density in {${\bf R}^k$}}, Israel J. Math., 54 (1986), pp.~307--316, \url{https://doi.org/10.1007/BF02764959}.

\bibitem{CastroSilva2021}
{\sc D.~Castro-Silva}, {\em Extremal sets with forbidden configurations and the independence ratio of geometric hypergraphs}, dissertation, University of Cologne, 2021, \url{https://kups.ub.uni-koeln.de/64364/}.

\bibitem{CastroSilva2023}
{\sc D.~Castro-Silva}, {\em Geometrical sets with forbidden configurations}, Forum Math. Sigma, 11 (2023), pp.~Paper No. e44, 45, \url{https://doi.org/10.1017/fms.2023.43}.

\bibitem{CastroSilvaOSV2022}
{\sc D.~Castro-Silva, F.~M. de~Oliveira~Filho, L.~Slot, and F.~Vallentin}, {\em A recursive {Lovász} theta number for simplex-avoiding sets}, Proc. Amer. Math. Soc., 150 (2022), pp.~3307--3322, \url{https://doi.org/10.1090/proc/15940}.

\bibitem{CastroSilvaOSV2023}
{\sc D.~Castro-Silva, F.~M. de~Oliveira~Filho, L.~Slot, and F.~Vallentin}, {\em A recursive theta body for hypergraphs}, Combinatorica, 43 (2023), pp.~909--938, \url{https://doi.org/10.1007/s00493-023-00040-9}.

\bibitem{CohenHLL1997}
{\sc G.~Cohen, I.~Honkala, S.~Litsyn, and A.~Lobstein}, {\em Covering codes}, vol.~54 of North-Holland Mathematical Library, North-Holland Publishing Co., Amsterdam, 1997.

\bibitem{CohnLL2024}
{\sc H.~Cohn, D.~de~Laat, and N.~Leijenhorst}, {\em Optimality of spherical codes via exact semidefinite programming bounds}, 2024, \url{https://arxiv.org/abs/2403.16874}.

\bibitem{CohnLS2022}
{\sc H.~Cohn, D.~de~Laat, and A.~Salmon}, {\em Three-point bounds for sphere packing}, 2022, \url{https://arxiv.org/abs/2206.15373}.

\bibitem{CohnE2003}
{\sc H.~Cohn and N.~Elkies}, {\em New upper bounds on sphere packings. {I}}, Ann. of Math. (2), 157 (2003), pp.~689--714, \url{https://doi.org/10.4007/annals.2003.157.689}.

\bibitem{CohnK2007}
{\sc H.~Cohn and A.~Kumar}, {\em Universally optimal distribution of points on spheres}, J. Amer. Math. Soc., 20 (2007), pp.~99--148, \url{https://doi.org/10.1090/S0894-0347-06-00546-7}.

\bibitem{CohnKMRV2017}
{\sc H.~Cohn, A.~Kumar, S.~D. Miller, D.~Radchenko, and M.~Viazovska}, {\em The sphere packing problem in dimension 24}, Ann. of Math. (2), 185 (2017), pp.~1017--1033, \url{https://doi.org/10.4007/annals.2017.185.3.8}.

\bibitem{CohnL2024}
{\sc H.~Cohn and A.~Li}, {\em Improved kissing numbers in seventeen through twenty-one dimensions}, 2024, \url{https://arxiv.org/abs/2411.04916}.

\bibitem{CohnS2021}
{\sc H.~Cohn and A.~Salmon}, {\em Sphere packing bounds via rescaling}, 2021, \url{https://arxiv.org/abs/2108.10936}.

\bibitem{ConwayS1999}
{\sc J.~H. Conway and N.~J.~A. Sloane}, {\em Sphere packings, lattices and groups}, vol.~290 of Grundlehren der mathematischen Wissenschaften [Fundamental Principles of Mathematical Sciences], Springer-Verlag, New York, third~ed., 1999, \url{https://doi.org/10.1007/978-1-4757-6568-7}.
\newblock With additional contributions by E. Bannai, R. E. Borcherds, J. Leech, S. P. Norton, A. M. Odlyzko, R. A. Parker, L. Queen and B. B. Venkov.

\bibitem{Croft1967}
{\sc H.~T. Croft}, {\em Incidence incidents}, Eureka, 30 (1967), pp.~22--26.

\bibitem{KlerkP2002}
{\sc E.~de~Klerk and D.~V. Pasechnik}, {\em Approximation of the stability number of a graph via copositive programming}, SIAM J. Optim., 12 (2002), pp.~875--892, \url{https://doi.org/10.1137/S1052623401383248}.

\bibitem{KlerkV2016}
{\sc E.~de~Klerk and F.~Vallentin}, {\em On the {Turing} model complexity of interior point methods for semidefinite programming}, SIAM J. Optim., 26 (2016), pp.~1944--1961, \url{https://doi.org/10.1137/15M103114X}.

\bibitem{Laat2020}
{\sc D.~de~Laat}, {\em Moment methods in energy minimization: new bounds for {Riesz} minimal energy problems}, Trans. Amer. Math. Soc., 373 (2020), pp.~1407--1453, \url{https://doi.org/10.1090/tran/7976}.

\bibitem{LaatOV2014}
{\sc D.~de~Laat, F.~M. de~Oliveira~Filho, and F.~Vallentin}, {\em Upper bounds for packings of spheres of several radii}, Forum Math. Sigma, 2 (2014), pp.~Paper No. e23, 42, \url{https://doi.org/10.1017/fms.2014.24}.

\bibitem{LaatLM2024}
{\sc D.~de~Laat, N.~M. Leijenhorst, and W.~H.~H. de~Muinck~Keizer}, {\em Optimality and uniqueness of the ${D}_4$ root system}, 2024, \url{https://arxiv.org/abs/2404.18794}.

\bibitem{LaatMM2022}
{\sc D.~de~Laat, F.~C. Machado, and W.~de~Muinck~Keizer}, {\em The {L}asserre hierarchy for equiangular lines with a fixed angle}, 2022, \url{https://arxiv.org/abs/2211.16471}.

\bibitem{LaatV2015}
{\sc D.~de~Laat and F.~Vallentin}, {\em A semidefinite programming hierarchy for packing problems in discrete geometry}, Math. Program., 151 (2015), pp.~529--553, \url{https://doi.org/10.1007/s10107-014-0843-4}.

\bibitem{OliveiraV2010}
{\sc F.~M. de~Oliveira~Filho and F.~Vallentin}, {\em Fourier analysis, linear programming, and densities of distance avoiding sets in {$\mathbb{R}^n$}}, J. Eur. Math. Soc. (JEMS), 12 (2010), pp.~1417--1428, \url{https://doi.org/10.4171/JEMS/236}.

\bibitem{OliveiraV2019}
{\sc F.~M. de~Oliveira~Filho and F.~Vallentin}, {\em A counterexample to a conjecture of {Larman} and {Rogers} on sets avoiding distance 1}, Mathematika, 65 (2019), pp.~785--787, \url{https://doi.org/10.1112/s0025579319000160}.

\bibitem{DeCorteOV2022}
{\sc E.~DeCorte, F.~M. de~Oliveira~Filho, and F.~Vallentin}, {\em Complete positivity and distance-avoiding sets}, Math. Program., 191 (2022), pp.~487--558, \url{https://doi.org/10.1007/s10107-020-01562-6}.

\bibitem{DeCorteP2016}
{\sc E.~DeCorte and O.~Pikhurko}, {\em Spherical sets avoiding a prescribed set of angles}, Int. Math. Res. Not. IMRN,  (2016), pp.~6095--6117, \url{https://doi.org/10.1093/imrn/rnv319}.

\bibitem{Delsarte1973}
{\sc P.~Delsarte}, {\em An algebraic approach to the association schemes of coding theory}, Philips Res. Rep. Suppl.,  (1973), pp.~vi+97.

\bibitem{DelsarteGS1977}
{\sc P.~Delsarte, J.~M. Goethals, and J.~J. Seidel}, {\em Spherical codes and designs}, Geometriae Dedicata, 6 (1977), pp.~363--388, \url{https://doi.org/10.1007/bf03187604}.

\bibitem{DezaL1997}
{\sc M.~M. Deza and M.~Laurent}, {\em Geometry of cuts and metrics}, vol.~15 of Algorithms and Combinatorics, Springer-Verlag, Berlin, 1997, \url{https://doi.org/10.1007/978-3-642-04295-9}.

\bibitem{DostertLM2021}
{\sc M.~Dostert, D.~de~Laat, and P.~Moustrou}, {\em Exact semidefinite programming bounds for packing problems}, SIAM J. Optim., 31 (2021), pp.~1433--1458, \url{https://doi.org/10.1137/20M1351692}.

\bibitem{Erdoes1985}
{\sc P.~Erdős}, {\em Problems and results in combinatorial geometry}, in Discrete geometry and convexity ({New York}, 1982), vol.~440 of Ann. New York Acad. Sci., New York Acad. Sci., New York, 1985, pp.~1--11, \url{https://doi.org/10.1111/j.1749-6632.1985.tb14533.x}.

\bibitem{FurstenbergKW1990}
{\sc H.~Furstenberg, Y.~Katznelson, and B.~Weiss}, {\em Ergodic theory and configurations in sets of positive density}, in Mathematics of {R}amsey theory, vol.~5 of Algorithms Combin., Springer, Berlin, 1990, pp.~184--198, \url{https://doi.org/10.1007/978-3-642-72905-8\_13}.

\bibitem{Ganzhinov2025}
{\sc M.~Ganzhinov}, {\em Highly symmetric lines}, Linear Algebra Appl., 722 (2025), pp.~12--37, \url{https://doi.org/10.1016/j.laa.2025.05.002}.

\bibitem{GrotschelLS1988}
{\sc M.~Grötschel, L.~Lovász, and A.~Schrijver}, {\em Geometric algorithms and combinatorial optimization}, vol.~2 of Algorithms and Combinatorics: Study and Research Texts, Springer-Verlag, Berlin, 1988, \url{https://doi.org/10.1007/978-3-642-97881-4}.

\bibitem{Hales2017}
{\sc T.~Hales, M.~Adams, G.~Bauer, T.~D. Dang, J.~Harrison, L.~T. Hoang, C.~Kaliszyk, V.~Magron, S.~McLaughlin, T.~T. Nguyen, Q.~T. Nguyen, T.~Nipkow, S.~Obua, J.~Pleso, J.~Rute, A.~Solovyev, T.~H.~A. Ta, N.~T. Tran, T.~D. Trieu, J.~Urban, K.~Vu, and R.~Zumkeller}, {\em A formal proof of the {Kepler} conjecture}, Forum Math. Pi, 5 (2017), pp.~e2, 29, \url{https://doi.org/10.1017/fmp.2017.1}.

\bibitem{Hales2005}
{\sc T.~C. Hales}, {\em A proof of the {Kepler} conjecture}, Ann. of Math. (2), 162 (2005), pp.~1065--1185, \url{https://doi.org/10.4007/annals.2005.162.1065}.

\bibitem{Halicka2002}
{\sc M.~Halická, E.~de~Klerk, and C.~Roos}, {\em On the convergence of the central path in semidefinite optimization}, SIAM J. Optim., 12 (2002), pp.~1090--1099, \url{https://doi.org/10.1137/S1052623401390793}.

\bibitem{HillW1987}
{\sc R.~D. Hill and S.~R. Waters}, {\em On the cone of positive semidefinite matrices}, Linear Algebra Appl., 90 (1987), pp.~81--88, \url{https://doi.org/10.1016/0024-3795(87)90307-7}.

\bibitem{Kalai2015}
{\sc G.~Kalai}, {\em Some old and new problems in combinatorial geometry {I}: around {Borsuk}'s problem}, in Surveys in combinatorics 2015, vol.~424 of London Math. Soc. Lecture Note Ser., Cambridge Univ. Press, Cambridge, 2015, pp.~147--174.

\bibitem{LarmanR1972}
{\sc D.~G. Larman and C.~A. Rogers}, {\em The realization of distances within sets in {Euclidean} space}, Mathematika, 19 (1972), pp.~1--24, \url{https://doi.org/10.1112/S0025579300004903}.

\bibitem{Lasserre2001}
{\sc J.~B. Lasserre}, {\em An explicit exact {SDP} relaxation for nonlinear 0-1 programs}, in Integer programming and combinatorial optimization ({Utrecht}, 2001), vol.~2081 of Lecture Notes in Comput. Sci., Springer, Berlin, 2001, pp.~293--303, \url{https://doi.org/10.1007/3-540-45535-3\_23}.

\bibitem{Laurent2003}
{\sc M.~Laurent}, {\em A comparison of the {Sherali-Adams}, {Lovász-Schrijver}, and {Lasserre} relaxations for 0-1 programming}, Math. Oper. Res., 28 (2003), pp.~470--496, \url{https://doi.org/10.1287/moor.28.3.470.16391}.

\bibitem{Laurent2007}
{\sc M.~Laurent}, {\em Strengthened semidefinite programming bounds for codes}, Math. Program., 109 (2007), pp.~239--261, \url{https://doi.org/10.1007/s10107-006-0030-3}.

\bibitem{LaurentV2026}
{\sc M.~Laurent and F.~Vallentin}, {\em A course on semidefinite optimization}, Cambridge University Press, Cambridge, in preparation.

\bibitem{LeijenhorstL2024}
{\sc N.~Leijenhorst and D.~de~Laat}, {\em Solving clustered low-rank semidefinite programs arising from polynomial optimization}, Math. Program. Comput., 16 (2024), pp.~503--534, \url{https://doi.org/10.1007/s12532-024-00264-w}.

\bibitem{LenstraLL}
{\sc A.~K. Lenstra, H.~W. Lenstra, Jr., and L.~Lovász}, {\em Factoring polynomials with rational coefficients}, Math. Ann., 261 (1982), pp.~515--534, \url{https://doi.org/10.1007/BF01457454}.

\bibitem{Levenshtein1979}
{\sc V.~I. Levenshtein}, {\em Boundaries for packings in {$n$}-dimensional {Euclidean} space}, Dokl. Akad. Nauk SSSR, 245 (1979), pp.~1299--1303.

\bibitem{Lovasz1979}
{\sc L.~Lovász}, {\em On the {Shannon} capacity of a graph}, IEEE Trans. Inform. Theory, 25 (1979), pp.~1--7, \url{https://doi.org/10.1109/TIT.1979.1055985}.

\bibitem{MachadoO2018}
{\sc F.~C. Machado and F.~M. de~Oliveira~Filho}, {\em Improving the semidefinite programming bound for the kissing number by exploiting polynomial symmetry}, Exp. Math., 27 (2018), pp.~362--369, \url{https://doi.org/10.1080/10586458.2017.1286273}.

\bibitem{McElieceRR1978}
{\sc R.~J. McEliece, E.~R. Rodemich, and H.~C. Rumsey, Jr.}, {\em The {Lovász} bound and some generalizations}, J. Combin. Inform. System Sci., 3 (1978), pp.~134--152.

\bibitem{MittelmannV2010}
{\sc H.~D. Mittelmann and F.~Vallentin}, {\em High-accuracy semidefinite programming bounds for kissing numbers}, Experiment. Math., 19 (2010), pp.~175--179, \url{https://doi.org/10.1080/10586458.2010.10129070}.

\bibitem{Monniaux2011}
{\sc D.~Monniaux and P.~Corbineau}, {\em On the generation of {P}ositivstellensatz witnesses in degenerate cases}, in Interactive theorem proving, vol.~6898 of Lecture Notes in Comput. Sci., Springer, Heidelberg, 2011, pp.~249--264, \url{https://doi.org/10.1007/978-3-642-22863-6\_19}.

\bibitem{Moser1991}
{\sc W.~O.~J. Moser}, {\em Problems, problems, problems}, vol.~31, 1991, pp.~201--225, \url{https://doi.org/10.1016/0166-218X(91)90071-4}.
\newblock First Canadian Conference on Computational Geometry (Montreal, PQ, 1989).

\bibitem{MotzkinS1965}
{\sc T.~S. Motzkin and E.~G. Straus}, {\em Maxima for graphs and a new proof of a theorem of {Turán}}, Canadian J. Math., 17 (1965), pp.~533--540, \url{https://doi.org/10.4153/CJM-1965-053-6}.

\bibitem{Musin2008}
{\sc O.~R. Musin}, {\em The kissing number in four dimensions}, Ann. of Math. (2), 168 (2008), pp.~1--32, \url{https://doi.org/10.4007/annals.2008.168.1}.

\bibitem{Nemirovski2007}
{\sc A.~Nemirovski}, {\em Advances in convex optimization: conic programming}, in International {Congress} of {Mathematicians}. {Vol. I}, Eur. Math. Soc., Zürich, 2007, pp.~413--444, \url{https://doi.org/10.4171/022-1/17}.

\bibitem{NieRS2010}
{\sc J.~Nie, K.~Ranestad, and B.~Sturmfels}, {\em The algebraic degree of semidefinite programming}, Math. Program., 122 (2010), pp.~379--405, \url{https://doi.org/10.1007/s10107-008-0253-6}.

\bibitem{Novikov2025}
{\sc A.~Novikov, N.~Vũ, M.~Eisenberger, E.~Dupont, P.-S. Huang, A.~Z. Wagner, S.~Shirobokov, B.~Kozlovskii, F.~J.~R. Ruiz, A.~Mehrabian, M.~P. Kumar, A.~See, S.~Chaudhuri, G.~Holland, A.~Davies, S.~Nowozin, P.~Kohli, and M.~Balog}, {\em Alphaevolve: A coding agent for scientific and algorithmic discovery}, 2025, \url{https://arxiv.org/abs/2506.13131}.

\bibitem{OdlyzkoS1979}
{\sc A.~M. Odlyzko and N.~J.~A. Sloane}, {\em New bounds on the number of unit spheres that can touch a unit sphere in {$n$} dimensions}, J. Combin. Theory Ser. A, 26 (1979), pp.~210--214, \url{https://doi.org/10.1016/0097-3165(79)90074-8}.

\bibitem{Putinar1993}
{\sc M.~Putinar}, {\em Positive polynomials on compact semi-algebraic sets}, Indiana Univ. Math. J., 42 (1993), pp.~969--984, \url{https://doi.org/10.1512/iumj.1993.42.42045}.

\bibitem{RienerRV2025}
{\sc C.~Riener, J.~Rolfes, and F.~Vallentin}, {\em A semidefinite programming hierarchy for covering problems in discrete geometry}, Numerical Algebra, Control and Optimization,  (2025), \url{https://doi.org/10.3934/naco.2025015}.

\bibitem{Schoenberg1942}
{\sc I.~J. Schoenberg}, {\em Positive definite functions on spheres}, Duke Math. J., 9 (1942), pp.~96--108, \url{http://projecteuclid.org/euclid.dmj/1077493072}.

\bibitem{Schrijver1979}
{\sc A.~Schrijver}, {\em A comparison of the {Delsarte} and {Lovász} bounds}, IEEE Trans. Inform. Theory, 25 (1979), pp.~425--429, \url{https://doi.org/10.1109/TIT.1979.1056072}.

\bibitem{Schrijver2005}
{\sc A.~Schrijver}, {\em New code upper bounds from the {Terwilliger} algebra and semidefinite programming}, IEEE Trans. Inform. Theory, 51 (2005), pp.~2859--2866, \url{https://doi.org/10.1109/TIT.2005.851748}.

\bibitem{SchutteW1953}
{\sc K.~Schütte and B.~L. van~der Waerden}, {\em Das {Problem} der dreizehn {Kugeln}}, Math. Ann., 125 (1953), pp.~325--334, \url{https://doi.org/10.1007/BF01343127}.

\bibitem{Simon2011}
{\sc B.~Simon}, {\em Convexity: An analytic viewpoint}, vol.~187 of Cambridge Tracts in Mathematics, Cambridge University Press, Cambridge, 2011, \url{https://doi.org/10.1017/CBO9780511910135}.

\bibitem{Szekely2002}
{\sc L.~A. Székely}, {\em Erdős on unit distances and the {Szemerédi-Trotter} theorems}, in Paul {Erdős} and his mathematics, {II} ({Budapest}, 1999), vol.~11 of Bolyai Soc. Math. Stud., János Bolyai Math. Soc., Budapest, 2002, pp.~649--666.

\bibitem{Tarski1948}
{\sc A.~Tarski}, {\em A {Decision} {Method} for {Elementary} {Algebra} and {Geometry}}, The Rand Corporation, Santa Monica, CA, 1948, \url{https://www.rand.org/content/dam/rand/pubs/reports/2008/R109.pdf}.

\bibitem{Viazovska2017}
{\sc M.~S. Viazovska}, {\em The sphere packing problem in dimension 8}, Ann. of Math. (2), 185 (2017), pp.~991--1015, \url{https://doi.org/10.4007/annals.2017.185.3.7}.

\bibitem{Witsenhausen1974}
{\sc H.~S. Witsenhausen}, {\em Spherical sets without orthogonal point pairs}, Amer. Math. Monthly, 81 (1974), pp.~1101--1102, \url{https://doi.org/10.2307/2319047}.

\bibitem{Zinoviev1999}
{\sc V.~A. Zinov’ev and T.~Ericson}, {\em New lower bounds for the kissing number in small dimensions}, Problemy Peredachi Informatsii, 35 (1999), pp.~3--11.

\end{thebibliography}
\end{document}